\documentclass[a4paper,twoside,10pt]{article}

\usepackage[colorlinks=true]{hyperref} 
\hypersetup{urlcolor=green,linkcolor=red,citecolor=blue,colorlinks=true}

\usepackage[frenchb]{babel}	
\usepackage[T1]{fontenc}	          
\usepackage[applemac]{inputenc}	 
\usepackage{vmargin}		 
\usepackage{fancyhdr}		 
\usepackage{calc}		           
\usepackage{ifthen}		           
\usepackage{pifont}		           
\usepackage{supertabular}	  
\usepackage{multicol}		  
\usepackage{wrapfig}		  
\usepackage{fancybox}		  
\usepackage{rotating}		  
\usepackage{xspace}		  

\usepackage{amsmath}	           
\usepackage{amssymb}	           
\usepackage{amsfonts}	           
\usepackage{verbatim}	           
\usepackage{mathrsfs}	           
\usepackage{dsfont}
\usepackage{amsthm}

\usepackage{array}
\usepackage{curves}
\usepackage{epic}
\usepackage{eepic}
\usepackage{epsfig}
\usepackage{graphicx}

\graphicspath{{PS/}}

\setpapersize{custom}{21.5cm}{26.5cm}
\setmarginsrb{43mm}{10mm}{43mm}{15mm}{10mm}{6mm}{0mm}{10mm}

\providecommand\StandardMathComma{\mathcode`\,="613B%
   \addto\extrasfrenchb{\mathcode`\,="613B}}

\StandardMathComma

\makeatletter
\declare@shorthand{frenchb}{:}{%
    \ifhmode
      \ifdim\lastskip>\z@
        \unskip\penalty\@M\thinspace
      \else
        \FDP@thinspace
      \fi
    \fi
    \string:}
\makeatother

\newcounter{ifnote}
\newcounter{IfAlignementEnonces}
\newcounter{Corrige}
\newcounter{Equation}[Corrige]

\newtheorem{theoreme}{Th\'eor\`eme}
\newtheorem*{theoreme imaginaire}{Th\'eor\`eme imaginaire}
\newtheorem{proposition}[theoreme]{Proposition}
\newtheorem{lemme}[theoreme]{Lemme}
\newtheorem{definition}[theoreme]{D\'efinition}
\newtheorem{corollaire}[theoreme]{Corollaire}

\newcommand{\demo}{
         \noindent
         {\it D\'emonstration :}
         }
\newcommand{\finpreuve}{\centersright[0]{}{$\Box$} \medskip}
\newcommand{\findemo}{\finpreuve}
\newcommand{\petiteremarque}{
         \noindent
         {\it Remarque :}
         }
\newcommand{\exemple}{
         \noindent
         {\it Exemple :}
         }

\newcommand{\exemplefondamental}{
         \noindent
         {\it Exemple fondamental :}
         }

\newcommand{\LeTitre}{}

\newcommand{\Titre}{
	\begin{center}
	{\bfseries\Large \LeTitre}
	\end{center}
	\vskip 1 cm
	}

\newcommand{\LeAuteur}{}

\newcommand{\titre}[1]{
	\renewcommand{\LeTitre}{#1}
	\Titre
	\stepcounter{Corrige}	
        \renewcommand{\labelitemi}{\ensuremath{\bullet}}
	}

\newlength{\PartieWidth}
\newcommand{\Partie}[3][8cm]{%
    \settowidth{\PartieWidth}{{\large\textbf{\scshape #3}}}%
    \ifthenelse{\lengthtest{\PartieWidth > #1}}{%
        \setlength{\PartieWidth}{#1}%
    }{}%
    \vspace*{\smallskipamount}%
    \begin{center}%
    {\large\textbf{#2.}}{\normalsize\quad}%
    \begin{minipage}[t]{\PartieWidth}%
    {\large\textbf{\scshape #3}}%
    \end{minipage}%
    \end{center}%
}

\newcounter{lareponse}

\newcommand{\enonce}[2][0]{%
        \noindent%
        \epsfig{file=#2,angle=#1,width=\linewidth,clip=}%
        \ifthenelse{\value{IfAlignementEnonces} = 1}{%
                \begin{picture}(0,0)%
                \multiput(-400,600)(0,-9){70}{\line(1,0){500}}%
                \end{picture}%
        }%
        {}%
}

\newlength{\Sujet}
\newlength{\Resume}
\newlength{\Enonces}
\newlength{\Corriges}
\newenvironment{sommaire}
	{\begin{supertabular}{
	        p{\Sujet}     
        	p{\Resume}       
	        p{\Enonces}      
		p{\Corriges}     
	        }
	}
	{\end{supertabular}}

\newcommand{\presentationCCP}[1]{
        \mbox{}
        \thispagestyle{empty}
        \vfill
        \centers{{\bfseries\Huge Concours Communs}}
        \centers{{\bfseries\Huge Polytechniques}}
        \mbox{}
        \vskip2cm
        \mbox{}
        \vfill
        \mbox{}
        \newpage
}

\newlength{\alignelength}

\mathcode`A="7041 \mathcode`B="7042 \mathcode`C="7043 \mathcode`D="7044
\mathcode`E="7045 \mathcode`F="7046 \mathcode`G="7047 \mathcode`H="7048
\mathcode`I="7049 \mathcode`J="704A \mathcode`K="704B \mathcode`L="704C
\mathcode`M="704D \mathcode`N="704E \mathcode`O="704F \mathcode`P="7050
\mathcode`Q="7051 \mathcode`R="7052 \mathcode`S="7053 \mathcode`T="7054
\mathcode`U="7055 \mathcode`V="7056 \mathcode`W="7057 \mathcode`X="7058
\mathcode`Y="7059 \mathcode`Z="705A

\catcode`\´=\active
\catcode`\ª=\active
\def´{\og\ignorespaces}
\defª{{\fg}}

        {\end{filet}\bigskip}

        {\end{filet}\bigskip}

\newenvironment{custom-itemize}[1]%
	{%
	\begin{list}{}%
		{%
		\settowidth{\labelwidth}{#1}%
		\setlength{\leftmargin}{\labelwidth+\labelsep}%
		}%
        }%
        {\end{list}}
	{\end{custom-itemize}}

\newcommand{\separation}{\begin{center}\rule{3 cm}{1 pt}\end{center}}

\reversemarginpar
\newcounter{Note}[page]
\newcounter{notesimple}
\newcounter{notemark}
\newcounter{notemarkref}
\newcounter{notetext}
\newcounter{notetextref}
\newcommand{\notemark}{%
	\ifthenelse{\equal{\value{notemark}}{\value{notemarkref}}}{%
        	\refstepcounter{notemark}%
        	\refstepcounter{notemarkref}%
        	\setcounter{Note}{\value{notemark}}%
       }{%
        	\refstepcounter{notemark}%
        	\setcounter{Note}{\value{notemark}}%
	}%
	\setcounter{notetextref}{\value{notemarkref}}%
	\begin{picture}(0,0)%
        	\put(-3,-3){\LARGE\theNote}
	\end{picture}\xspace
}

\newcommand{\notetextbase}[1]{%
        \ifthenelse{\value{ifnote} = 1}{%
        \marginpar{\notedebasenumero{#1}}%
        }{}%
}

\newcommand{\notetext}[1]{
	\ifthenelse{\equal{\value{notetext}}{\value{notetextref}}}{%
        	\refstepcounter{notetext}%
        	\refstepcounter{notetextref}%
        	\setcounter{Note}{\value{notetext}}
        }{
        	\refstepcounter{notetext}%
        	\setcounter{Note}{\value{notetext}}
	}%
\setcounter{notemarkref}{\value{notetextref}}%
\notetextbase{#1}%
}

\newcommand{\notedebase}[1]{%
	\raggedright%
        \footnotesize%
        \vskip-0.5\baselineskip
        \rule[-1.4mm]{\linewidth}{0.5pt}
        \rule[1.4mm]{\linewidth}{0.5pt} \\%
        \vskip-0.5\baselineskip
        #1 \\%
        \vskip-\baselineskip
        \rule[-1.4mm]{\linewidth}{0.5pt}
        \rule[1.4mm]{\linewidth}{0.5pt} \\%
}

\newcommand{\notedebasenumero}[1]{%
	{\LARGE\theNote}
	\notedebase{#1}%
}

\newcommand{\notesimple}[1]{%
	\ifthenelse{\value{ifnote} = 1}{%
	\marginpar{\notedebase{#1}}%
	}{}%
}

\newcommand{\notenumero}[1]{%
        \ifthenelse{\value{ifnote} = 1}{%
        \refstepcounter{Note}%
        \setcounter{notemarkref}{\value{Note}}%
        \setcounter{notetextref}{\value{Note}}%
        \setcounter{notemark}{\value{Note}}%
        \setcounter{notetext}{\value{Note}}%
        \begin{picture}(0,0)%
        \put(-3,-3){\LARGE\theNote}
        \end{picture}%
        \ifthenelse{\isodd{\value{Note}}}{%
        	\protect\reversemarginpar%
	      	\marginpar[{\notedebasenumero{#1}}]{\notedebasenumero{#1}}}{%
           	\protect\normalmarginpar%
       		\marginpar[{\notedebasenumero{#1}}]{\notedebasenumero{#1}}}%
        }{}%
}

\newcommand{\note}[2][1]{%
	\ifthenelse{\equal{#1}{0} \or \equal{#1}{simple}}{%
                \setcounter{notesimple}{1}%
        }{}%
	\ifthenelse{\value{notesimple} = 0}{%
		\notenumero{#2}%
	}{%
		\notesimple{#2}%
	}
}

\newcounter{cinquantaine}
\newcounter{centaine}
\newcommand{\reperes}[1][20]{
	\ifthenelse{\value{ifnote} = 1}{%
		\setcounter{cinquantaine}{#1/5+1}%
		\setcounter{centaine}{#1/10+1}%
		\multiput(0,0)(10,0){#1}{\line(0,1){1}}%
		\multiput(0,0)(10,0){#1}{\line(0,-1){1}}%
		\multiput(0,0)(0,10){#1}{\line(1,0){1}}%
		\multiput(0,0)(0,10){#1}{\line(-1,0){1}}%
		\multiput(0,0)(-10,0){#1}{\line(0,1){1}}%
		\multiput(0,0)(-10,0){#1}{\line(0,-1){1}}%
		\multiput(0,0)(0,-10){#1}{\line(1,0){1}}%
		\multiput(0,0)(0,-10){#1}{\line(-1,0){1}}%
		\multiput(0,0)(50,0){\thecinquantaine}{\line(0,1){3}}%
		\multiput(0,0)(50,0){\thecinquantaine}{\line(0,-1){3}}%
		\multiput(0,0)(0,50){\thecinquantaine}{\line(1,0){3}}%
		\multiput(0,0)(0,50){\thecinquantaine}{\line(-1,0){3}}%
		\multiput(0,0)(-50,0){\thecinquantaine}{\line(0,1){3}}%
		\multiput(0,0)(-50,0){\thecinquantaine}{\line(0,-1){3}}%
		\multiput(0,0)(0,-50){\thecinquantaine}{\line(1,0){3}}%
		\multiput(0,0)(0,-50){\thecinquantaine}{\line(-1,0){3}}%
		\multiput(0,0)(100,0){\thecentaine}{\line(0,1){5}}%
		\multiput(0,0)(100,0){\thecentaine}{\line(0,-1){5}}%
		\multiput(0,0)(0,100){\thecentaine}{\line(1,0){5}}%
		\multiput(0,0)(0,100){\thecentaine}{\line(-1,0){5}}%
		\multiput(0,0)(-100,0){\thecentaine}{\line(0,1){5}}%
		\multiput(0,0)(-100,0){\thecentaine}{\line(0,-1){5}}%
		\multiput(0,0)(0,-100){\thecentaine}{\line(1,0){5}}%
		\multiput(0,0)(0,-100){\thecentaine}{\line(-1,0){5}}%
	}%
	{}%
}

\newcommand{\mathgras}[1]{\ensuremath{%
	{\text{\mathversion{bold}\ensuremath{#1}}}%
	}}

\renewcommand{\thesection}{\Roman{section}}

\newcommand{\seclabel}[1]{{%
  \renewcommand{\thesection}{\thechapter.\Roman{section}}%
  \addtocounter{section}{-1}\refstepcounter{section}%
  \label{#1}}}

\def\leaderfill{\leaders\hbox to 1ex{\hss.\hss}\hfill}

        {\setlength{\columnsep}{8 mm}%
        \setlength{\columnseprule}{0.5 pt}%
        \begin{multicols}{2}}%
        {\end{multicols}}

\makeatletter

\newbox\bk@bxb
\newbox\bk@bxa
\newif\if@bkcont 
\newif\ifbkcount
\newcount\bk@lcnt

\def\breakboxskip{2pt}
\def\breakboxparindent{1.8em}
\def\margesep{1cm}	
\def\intervalle{1mm}	

\def\filet{\vskip\breakboxskip\relax
\setbox\bk@bxb\vbox\bgroup
\advance\linewidth -\fboxrule
\advance\linewidth -\margesep
\advance\linewidth -\intervalle
\advance\linewidth -\fboxsep
\hsize\linewidth\@parboxrestore
\parindent\breakboxparindent\relax}

\def\bk@split{%
\@tempdimb\ht\bk@bxb 
\advance\@tempdimb\dp\bk@bxb 
\setbox\bk@bxa\vsplit\bk@bxb to\z@ 
\setbox\bk@bxa\vbox{\unvbox\bk@bxa}
\setbox\@tempboxa\vbox{\copy\bk@bxa\copy\bk@bxb}
\advance\@tempdimb-\ht\@tempboxa 
\advance\@tempdimb-\dp\@tempboxa}

\def\bk@addfsepht{%
     \setbox\bk@bxa\vbox{\vskip\fboxsep\box\bk@bxa}}

\def\bk@addskipht{%
     \setbox\bk@bxa\vbox{\vskip\@tempdimb\box\bk@bxa}}

\def\bk@addfsepdp{%
     \@tempdima\dp\bk@bxa
     \dp\bk@bxa\@tempdima}

\def\bk@addskipdp{%
     \@tempdima\dp\bk@bxa
     \advance\@tempdima\@tempdimb
     \dp\bk@bxa\@tempdima}

\def\bk@line{%
    \hbox to \linewidth{\ifbkcount\smash{\llap{\the\bk@lcnt\ }}\fi
    \hskip\margesep
    \vrule \@width\fboxrule\hskip\fboxsep
    \hskip\intervalle
    \box\bk@bxa\hfil
	}}

\def\endfilet{\egroup
\ifhmode\par\fi{\noindent\bk@lcnt\@ne 
\@bkconttrue\baselineskip\z@\lineskiplimit\z@
\lineskip\z@\vfuzz\maxdimen
\bk@split\bk@addfsepht\bk@addskipdp
\ifvoid\bk@bxb      
\def\bk@fstln{\bk@addfsepdp
\vbox{\bk@line}}%
\else               
\def\bk@fstln{\vbox{\bk@line}\hfil
\advance\bk@lcnt\@ne
\loop 
 \bk@split\bk@addskipdp\leavevmode
\ifvoid\bk@bxb      
 \@bkcontfalse\bk@addfsepdp
 \vtop{\bk@line}%
\else               
 \bk@line
\fi
 \hfil\advance\bk@lcnt\@ne
\if@bkcont\repeat}%
\fi
\leavevmode\bk@fstln\par}\vskip\breakboxskip\relax}

\bkcountfalse

\makeatother

\newlength{\leftlength}
\newlength{\rightlength}
\newlength{\calculskip}
\newcommand{\calculvskip}[1]{%
  \ifthenelse{#1 = 0}{\setlength{\calculskip}{0pt}}{}%
  \ifthenelse{#1 = 1}{\setlength{\calculskip}{\smallskipamount}}{}%
  \ifthenelse{#1 = 2}{\setlength{\calculskip}{\medskipamount}}{}%
  \ifthenelse{#1 = 3}{\setlength{\calculskip}{\bigskipamount}}{}%
  \ifthenelse{#1 = 4}{\setlength{\calculskip}{1cm}}{}%
  \vskip\calculskip
}

\newcommand{\leftcentersright}[4][2]{%
        \settowidth{\leftlength}{#2}%
        \settowidth{\rightlength}{#4}%
        \calculvskip{#1}
	\noindent#2\hskip-\leftlength%
	\hskip\linewidth\hskip-\textwidth%
	\hfill#3\hfill
	\mbox{}\hskip-\rightlength#4%
        \vskip\calculskip%
	}

\newcommand{\centers}[2][2]{\leftcentersright[#1]{}{#2}{}}

\newcommand{\centersright}[3][2]{\leftcentersright[#1]{}{#2}{#3}}

\newsavebox{\boite}
\def\debutcom{\begin{lrbox}{\boite}}
\def\fincomg{\end{lrbox}\makebox[0cm][l]{\usebox{\boite}}%
             \hskip\linewidth\hskip-\textwidth}
\def\fincomd{\end{lrbox}\makebox[0cm][r]{\usebox{\boite}}}

\newenvironment{calculs:base}[2][2]{%
  	\calculvskip{#1}
  	\noindent
  	\begin{tabular*}{\linewidth}[t]%
    		{@{}>{\debutcom}l<{\fincomg}@{\extracolsep{\fill}}%
      		>{$}r<{$}%
      		@{$\ #2\ $}%
      		@{\extracolsep{0pt}}>{$}l<{$}%
      		@{\extracolsep{\fill}}>{\debutcom}r<{\fincomd}@{}}
	}{%
  	\end{tabular*}%
  	\vskip\calculskip
}

\newenvironment{inegalites:leq}[1][2]{%
	\begin{calculs:base}[#1]{\leq}}{%
	\end{calculs:base}
}

\newenvironment{inegalites:geq}[1][2]{%
	\begin{calculs:base}[#1]{\geq}}{%
	\end{calculs:base}
}

        {\begin{calculs:rcl}[#1]{r}{c}{l}}%
        {\end{calculs:rcl}}

\newenvironment{calculs:rcl}[4][2]{%
	\calculvskip{#1}
  	\noindent
  	\begin{tabular*}{\linewidth}[t]%
  		{@{}>{\debutcom}l<{\fincomg}@{\extracolsep{\fill}}
		>{$}#2<{$}@{\extracolsep{0pt}}%
		>{$\ }#3<{\ $}%
      		@{\extracolsep{0pt}}>{$}#4<{$}%
      		@{\extracolsep{\fill}}>{\debutcom}r<{\fincomd}@{}}%
	}{%
  	\end{tabular*}%
  	\vskip\calculskip
}

\newenvironment{calculs:rcl:extracol}[5][2]{%
	\calculvskip{#1}
 	\noindent
	\begin{tabular*}{\linewidth}[t]%
                {@{}>{\debutcom}l<{\fincomg}@{\extracolsep{\fill}}%
		>{$}#2<{$}@{\extracolsep{0pt}}%
		>{$\ }#3<{\ $}@{\extracolsep{0pt}}%
		>{$}#4<{$}@{\extracolsep{0pt}}%
      		>{$\null}#5<{$}%
      		@{\extracolsep{\fill}}>{\debutcom}r<{\fincomd}@{}}%
}{%
  	\end{tabular*}%
  	\vskip\calculskip
}

\newenvironment{calculs:latotale}[6][2]{%
  	\calculvskip{#1}
  	\noindent
  	\begin{tabular*}{\linewidth}[t]%
    	        {@{}>{\debutcom}l<{\fincomg}@{\extracolsep{\fill}}%
      		>{$\null}#2<{$}@{\extracolsep{0pt}}%
      		>{$\null}#3<{$}@{\extracolsep{0pt}}%
      		>{$\null}#4<{$}%
      		@{\extracolsep{0pt}}>{$\null}#5<{$}%
     		@{\extracolsep{0pt}}>{$\null}#6<{$}%
      		@{\extracolsep{\fill}}>{\debutcom}r<{\fincomd}@{}}%
}{%
  	\end{tabular*}%
  	\vskip\calculskip
}

\newcommand{\tf}[2]{\ensuremath{#1/#2}}
\newcommand{\tfi}[2]{\ensuremath{#1\textbackslash#2}}

\newcommand{\pa}[1]{\ensuremath{\left(#1\right)}}

\newcommand{\crochets}[1]{\ensuremath{\left[#1\right]}}
\newcommand{\pac}[1]{\ensuremath{\crochets{#1}}}		

\newcommand{\accolades}[1]{\ensuremath{\left\{#1\right\}}}
\newcommand{\paa}[1]{\ensuremath{\accolades{#1}}}		

\renewcommand{\tilde}[1]{\ensuremath{\widetilde{#1}}}

\newcommand{\plus}{\mbox{\protect\raisebox{.2mm}{\tiny{\ensuremath{+}}}}}
\newcommand{\moins}{\mbox{\protect\raisebox{.2mm}{\tiny{\ensuremath{-}}}}}
\newcommand{\pinf}{\plus\ensuremath{\infty}}
\newcommand{\minf}{\moins\ensuremath{\infty}}

\newlength{\pmlength} 
\newlength{\pmgraslength} 
\renewcommand{\leq}{\ensuremath{\leqslant}}
\renewcommand{\geq}{\ensuremath{\geqslant}}
\renewcommand{\epsilon}{\ensuremath{\varepsilon}}

\renewcommand{\setminus}{\smallsetminus}	

\newcommand{\qetq}{\quad\text{et}\quad}
\newcommand{\qqetqq}{\qquad\text{et}\qquad}

\newcommand{\qsiq}{\quad\text{si}\quad}

\newcommand{\fonction}[5]{%
        \ensuremath{#1\colon
        \left\{\hskip -1.5 mm
        \begin{array}{c@{\ }c@{\ }l}
        \medskip #2 & \longrightarrow & #3 \\
        #4 & \longmapsto & #5 \\
        \end{array}
        \right .
        }}

\newcommand{\application}[4]{
        \begin{array}{c@{\ }c@{\ }l}
        \medskip #1 & \longrightarrow & #2 \\
        #3 & \longmapsto & #4 \\
        \end{array}
        }

\newlength{\restsubwidth}
\newlength{\restsubheight}
\newlength{\restsubmoreheight}
\newcommand{\rest}[2]{%
        \settowidth{\restsubwidth}{\ensuremath{#2}}
        \settoheight{\restsubheight}{\ensuremath{{}_{#2}}}
        \ensuremath{{#1\hskip 0.5 pt}_{\vrule\kern2pt\parbox[b][%
        \the\restsubheight +
                \the\restsubmoreheight][b]{\the\restsubwidth}{%
                        \ensuremath{{}_{#2}}}}}
        }

\newcommand{\Sumtproto}[2]{%
        \ifthenelse{%
                \equal{#1}{}
        }{%
                \ifthenelse{%
                        \equal{#2}{}%
                }{%
                        \ensuremath{\sum}%
                }{%
                        \smash[b]{\ensuremath{\sum\limits_{#1}^{#2}}}%
                }
        }{%
                \ensuremath{\sum\limits}_{#1}^{#2}%
        }%
}

\newcommand{\union}[2]{\ensuremath{\bigcup\limits_{#1}^{#2}}}

\makeatletter
\newif\if@ListeStar
\global\@ListeStartrue

\newcommand{\liste}{%
	\@ifstar{\global\@ListeStartrue\@liste}%
		{\global\@ListeStarfalse\@liste}%
}

\newcommand{\@liste}[2][n]{%
	\if@ListeStar%
		\left({#2}_0,{#2}_1,\ldots,{#2}_{#1}\right)%
	\else%
		\left({#2}_1,{#2}_2,\ldots,{#2}_{#1}\right)%
	\fi\@ListeStarfalse%
}
\makeatother

\makeatletter
\newif\if@SuiteStar
\global\@SuiteStartrue

\newcommand{\suite}{%
	\@ifstar{\global\@SuiteStartrue\@suite}%
		{\global\@SuiteStarfalse\@suite}%
}

\newcommand{\@suite}[2][n]{%
	\if@SuiteStar%
		\left(#2_{#1}\right)_{#1\in\N^*}%
	\else%
		\left(#2_{#1}\right)_{#1\in\N}%
	\fi\@SuiteStarfalse%
}
\makeatother

\makeatletter
\newif\if@laststared
\global\@laststaredtrue
\newcommand{\mathBB}[1]{%
        \@ifstar%
        {\global\@laststaredtrue\m@thBB{#1}}%
        {\global\@laststaredfalse\m@thBB{#1}}%
}
\newcommand{\m@thBB}[1]{%
	\if@laststared{\ensuremath{\mathbb{#1}^{*}}\xspace}%
	\else{\ensuremath{\mathbb{#1}}\xspace}%
	\fi%
	\@laststaredfalse%
}
\makeatother

\renewcommand{\emptyset}{\ensuremath{\varnothing}}	
\newcommand{\N}{\ensuremath{\mathBB{N}}}		         
\newcommand{\Z}{\ensuremath{\mathBB{Z}}}		         
\newcommand{\Nstar}{\ensuremath{\N*}\xspace}

\newcommand{\TestGauche}[1]{\ifthenelse{\equal{#1}{}}{\minf}{#1}}
\newcommand{\TestDroite}[1]{\ifthenelse{\equal{#1}{}}{\pinf}{#1}}

\newcommand{\intn}[2]{\ensuremath{[\![ \, #1 \,;\, #2 \,]\!]}}
\newcommand{\eq}{\ensuremath{\Longleftrightarrow}}
	{\begin{pmatrix}}%
	{\end{pmatrix}}
	{\begin{vmatrix}}%
	{\end{vmatrix}}

\newlength{\moreinterligne}
\newlength{\letterarrowvskip}
\newlength{\vectlength}
\newlength{\vectheight}
\newlength{\boxlength}
\newcommand{\vect}[1]{%
        \settowidth{\vectlength}{\ensuremath{#1}}%
        \settoheight{\vectheight}{\ensuremath{#1}}%
        \settowidth{\boxlength}{%
	        \ensuremath{%
        	\overrightarrow{%
		        \parbox[b][\the\vectheight + \letterarrowvskip][b]{%
				\the\vectlength%
			}{\ensuremath{#1}}%
		}%
		}%
        }
        \parbox[b][\the\vectheight + \moreinterligne][b]{\the\boxlength}{%
	        \ensuremath{%
		        \overset{\hbox to \the\boxlength{\rightarrowfill}}{%
	        	\parbox[b][\the\vectheight + \letterarrowvskip][b]{
				\the\vectlength%
			}{\ensuremath{#1}}%
			}%
		}%
        }
}

\newlength{\smoreinterligne}
\newlength{\sletterarrowvskip}
\newlength{\svectlength}
\newlength{\svectheight}
\newlength{\sboxlength}

\newcommand{\svect}[1]{%
        \settowidth{\svectlength}{\mbox{\scriptsize{\ensuremath{#1}}}}%
        \settoheight{\svectheight}{\mbox{\scriptsize{\ensuremath{#1}}}}%
        \settowidth{\sboxlength}{%
	        \ensuremath{%
	        \overrightarrow{%
		        \parbox[b][\the\svectheight + \sletterarrowvskip][b]{%
				\the\svectlength%
			}%
        		{\mbox{\scriptsize{\ensuremath{#1}}}}%
		}%
		}%
        }
	\parbox[b][\the\svectheight + \smoreinterligne][b]{\the\sboxlength}{%
	        \ensuremath{%
		\overrightarrow{%
		        \parbox[b][\the\svectheight + \sletterarrowvskip][b]{%
				\the\svectlength%
			}%
	        	{\mbox{\scriptsize{\ensuremath{#1}}}}%
		}%
		}%
	}
}

	{\ensuremath{\left \{ \hskip -1.5 mm \begin{array}{#1@{\ =\ }l}}}%
	{\end{array}\right.}
	{\ensuremath{\begin{array}[t]{#1@{\ =\ }l}}}%
	{\end{array}}
	{\ensuremath{\left \{ \hskip -1.5 mm%
	 \begin{array}{#1@{\ }c@{\ }l}}}%
	{\end{array}\right.}
	{\ensuremath{\left \{ \hskip -1.5 mm \begin{array}{#1@{\quad}l}}}%
	{\end{array}\right.}
	{\ensuremath{\begin{array}[t]{#1@{\ \leq\ }l}}}%
	{\end{array}}
	{\ensuremath{\begin{array}[t]{#1@{\ }c@{\ }l}}}%
	{\end{array}}

\newlength{\boxrulewidth}
\newlength{\boxrulesep}
\newlength{\gauchelong}
\newlength{\droitelong}
\newlength{\colonnealong}
\newlength{\colonneblong}
\newlength{\colonneclong}

\newlength{\maxhaut}
\newlength{\maxbas}

\newcommand{\writeifexist}[1]{%
        \ifthenelse{\equal{#1}{\null}}{\null}{#1&}%
}

\newcommand{\lengthifexist}[2]{%
       \ifthenelse{\not\equal{#1}{\null}}{%
                \settowidth{#2}{$#1$}
                }{\setlength{#2}{0cm}}%
}

\makeatletter

\def\arraybox#1{
    \@ifnextchar[
    	{\iarraybox{#1}}%
    	{\ivarraybox{#1}{\null}{\null}[\null]}}

\def\iarraybox#1[#2]{%
    \@ifnextchar[
    	{\iiarraybox{#1}{#2}}%
    	{\ivarraybox{#1}{#2}{\null}[\null]}}

\def\iiarraybox#1#2[#3]{%
    \@ifnextchar[
  	  {\iiiarraybox{#1}{#2}{#3}}%
    	{\ivarraybox{#1}{#2}{#3}[\null]}}

\def\iiiarraybox#1#2#3[#4]{%
    \@ifnextchar[
  	  {\ivarraybox{#1}{#2}{#3}[#4]\relax}%
    	{\ivarraybox{#1}{#2}{#3}[#4]}}%

\def\ivarraybox#1#2#3[#4]#5{%
        \settowidth{\gauchelong}{$#1\ $}%
        \lengthifexist{#2}{\colonnealong}%
        \lengthifexist{#3}{\colonneblong}%
        \lengthifexist{#4}{\colonneclong}%
        \settoheight{\maxhaut}{$#1#2#3#4#5$}
        \settodepth{\maxbas}{$#1#2#3#4#5$}
        \raisebox{-\maxbas-\boxrulewidth-\boxrulesep}{%
                \rule{\boxrulewidth}{%
                \maxbas+\maxhaut+2\boxrulewidth+2\boxrulesep}}%
        \raisebox{\maxhaut+\boxrulesep}{%
                \makebox[0cm][l]{\rule{%
                \gauchelong+\colonnealong+\colonneblong+\colonneclong}{%
                \boxrulewidth}}}%
        \raisebox{-\maxbas-\boxrulewidth-\boxrulesep}{%
                \makebox[0cm][l]{\rule{%
                \gauchelong+\colonnealong+\colonneblong+\colonneclong}{%
                \boxrulewidth}}}%
        \hskip\boxrulesep%
        #1&%
        \writeifexist{#2}%
        \writeifexist{#3}%
        \writeifexist{#4}%
        #5%
        \settowidth{\droitelong}{$#5\ \eq$}
        \settoheight{\maxhaut}{$#1#2#3#4#5$}
        \settodepth{\maxbas}{$#1#2#3#4#5$}
        \hskip\boxrulesep%
        \hskip-\droitelong%
        \raisebox{\maxhaut+\boxrulesep}{%
                \makebox[0cm][l]{\rule{\droitelong}{\boxrulewidth}}}
        \raisebox{-\maxbas-\boxrulesep-\boxrulewidth}{%
                \makebox[0cm][l]{\rule{\droitelong}{\boxrulewidth}}}
        \hskip\droitelong%
        \raisebox{-\maxbas-\boxrulesep-\boxrulewidth}{%
                \rule{\boxrulewidth}{%
                \maxbas+\maxhaut+2\boxrulesep+2\boxrulewidth}}
}

\makeatother

\newlength{\angstromlength}
\usepackage[all]{xy}

\begin{document}
\titre{Th\'eor\`eme de Kurosh pour les relations d'\'equivalence bor\'eliennes}
\centers{\textbf{Aur\'elien Alvarez}}
\centers{\texttt{aurelien.alvarez@epfl.ch}}

\bigskip
\bigskip

En th\'eorie des groupes, le th\'eor\`eme de Kurosh (\cite{MR0071422}) est un r\'esultat de structure concernant les sous-groupes d'un produit libre de groupes ; plus pr\'ecis\'ement, un sous-groupe $H$ du produit libre~$G$ d'une famille de groupes $(G_z)_{z \in Z}$ est isomorphe au produit libre de son intersection avec des conjugu\'es des $G_z$ convenablement index\'es et d'un sous-groupe libre de $G$. Le th\'eor\`eme principal de cet article (th. \ref{theoreme de Kurosh}) est un r\'esultat analogue dans le cadre des relations d'\'equivalence bor\'eliennes \`a classes d\'enombrables, que nous d\'emontrons en d\'eveloppant une th\'eorie de Bass-Serre dans ce cadre particulier. Nous renvoyons \`a \cite{Alv08b} pour une th\'eorie de Bass-Serre dans le cadre plus naturel des groupo\"{i}des bor\'eliens.

\separation

\'Etant donn\'e une relation d'\'equivalence bor\'elienne $\cal R$ \`a classes d\'enombrables sur un espace bor\'elien standard $X$, les acteurs principaux de ce travail sont les $\cal R$-{\it arboretums} (d\'ef. \ref{definition arboretum}), c'est-\`a-dire la donn\'ee d'une {\it action} de $\cal R$ sur un {\it champ d'arbres bor\'elien} sur $X$. Nous nous int\'eressons dans un premier temps aux actions {\it quasi-libres}, ce qui nous permet d'obtenir une d\'emonstration g\'eom\'etrique qu'une sous-relation d'une relation d'\'equivalence bor\'elienne arborable est arborable (cor. \ref{sous-relation d'une arborable}, voir aussi \cite{MR1900547} dans le cadre bor\'elien et \cite{MR1728876} en pr\'esence d'une mesure).

\medskip

\`A toute d\'ecomposition de $\cal R$ en produit libre de deux sous-relations ${\cal R}_1$ et ${\cal R}_2$, est canoniquement associ\'e un $\cal R$-arboretum {\it bi-color\'e} et le th\'eor\`eme \ref{caract\'erisation des produits amalgam\'es} donne une caract\'erisation {\it dynamique} des produits amalgam\'es de deux sous-relations suivant une sous-relation commune. Via la notion de {\it graphe de relations} (d\'ef. \ref{definition graphe de relations}), nous d\'emontrons l'existence d'une {\it d\'esingularisation} pour toute action de $\cal R$ sur un arboretum (th. \ref{desingularisation d'une action}) et donnons des r\'esultats sur la structure de $\cal R$ (prop. \ref{produits amalgames}, \ref{arborable} et~\ref{generateurs}). Enfin, nous utilisons les r\'esultats obtenus pour d\'emontrer un analogue du th\'eor\`eme de Kurosh pour les sous-relations d'une relation d'\'equivalence bor\'elienne ${\cal R} = {\star}_{i \in I} {\cal R}_i$ qui se d\'ecompose en produit libre d\'enombrable de sous-relations ${\cal R}_i$.

\begin{theoreme}[th. \ref{theoreme de Kurosh}]
Soit $\cal R$ une relation d'\'equivalence bor\'elienne sur $X$, produit libre d\'enombrable de sous-relations ${\cal R}_i$ ($i \in I$). Alors
\centers{${\cal S} = {\star}_{i \in I}\pa{{\star}_{k_i \in K(i)} {\cal S}_{k_i}} \star {\cal T},$}
\noindent
o\`u, pour tout $k_i$ d'un ensemble d\'enombrable $K(i)$, il existe un \'el\'ement $\phi_{k_i}$ de $[[\cal R]]$ d\'efini sur une partie bor\'elienne $A_{k_i}$ de $X$ tel que
\centers{${\cal S}_{k_i} = \pa{\phi_{k_i}^{-1} \rest{{\cal R}_{i}}{\phi_{k_i}(A_{k_i})} \phi_{k_i}} \cap {\cal S},$}
\noindent
et o\`u $\cal T$ est une sous-relation arborable de $\cal S$. De plus, pour tout $i$ de $I$, il existe $k_i$ dans~$K(i)$ tels que
\centers{$A_{k_i}=X \qqetqq {\cal S}_{k_i}={\cal R}_i \cap {\cal S}.$}
\end{theoreme}

En particulier, nous donnons la d\'ecomposition de la restriction de $\cal R$ \`a toute partie bor\'elienne~$Y$ de $X$ (th. \ref{theoreme restriction}) et pr\'ecisons ainsi les r\'esultats de Ioana-Peterson-Popa (\cite{ioana-peterson-popa}) obtenus dans le cas de facteurs ergodiques.

\bigskip

En collaboration avec D. Gaboriau, nous introduisons dans \cite{MFI} la notion de relation d'\'equivalence mesur\'ee {\it librement ind\'ecomposable} ainsi que la classe des groupes d\'enombrables {\it mesurablement librement ind\'ecomposables}. Les d\'emonstrations des r\'esultats de rigidit\'e que nous obtenons (th. IV.18 et th. V.1) reposent en grande partie sur les th\'eor\`emes \ref{theoreme de Kurosh} et \ref{theoreme restriction} de cet article.

\bigskip

\noindent
{\it Remerciements : }
Je tiens \`a remercier sinc\`erement Damien Gaboriau pour son encouragement tout au long de ce travail ainsi que Fr\'ed\'eric Paulin pour tout le soin qu'il a accord\'e \`a une premi\`ere version de ce texte.

\separation

Dans la suite, le couple $(X,{\cal B}_X)$ d\'esigne toujours un espace bor\'elien standard et $\cal R$ une relation d'\'equivalence bor\'elienne \`a classes d\'enombrables sur $X$.

\section{Actions quasi-libres et arboralit\'e}\label{paragraphe quasi-libre}

Nous commen\c{c}ons par rappeler la d\'efinition des $\cal R$-{\it espaces fibr\'es standards} et mentionnons quelques propri\'et\'es de ces derniers qui nous seront utiles par la suite. Nous introduisons ensuite les $\cal R$-arboretums (d\'ef. \ref{definition arboretum}) et nous nous int\'eressons au cas particulier d'actions quasi-libres (d\'ef.~\ref{definition quasi-libre}). Nous obtenons ainsi une caract\'erisation dynamique des relations d'\'equivalence bor\'eliennes arborables (th. \ref{quasi-libre=arborable}).

\smallskip

Rappelons qu'une partie bor\'elienne $A$ de $X$ est un {\it domaine fondamental} de $\cal R$ si elle rencontre chaque classe de $\cal R$ en un unique \'el\'ement et qu'une relation d'\'equivalence bor\'elienne est lisse si elle admet un domaine fondamental. Le {\it satur\'e} ${\cal R} \cdot A$ d'une partie bor\'elienne $A$ de $X$ est la partie bor\'elienne de $X$ constitu\'ee des \'el\'ements $\cal R$-\'equivalents \`a un \'el\'ement de $A$. Lorsque le satur\'e de $A$ co\"{i}ncide avec~$X$, on dit que $A$ est un {\it domaine complet} de $\cal R$. Deux relations d\'equivalence bor\'eliennes~$\cal R$ et~$\cal R'$ sur $X$ et $X'$ respectivement sont {\it stablement orbitalement \'equivalentes} (ou {\it stablement isomorphes}) s'il existe des domaines complets $A$ de $\cal R$ et~$A'$ de~${\cal R}'$ telles que les restrictions de $\cal R$ et de ${\cal R}'$ \`a ces domaines complets soient orbitalement \'equivalentes.  Le {\it pseudo-groupe plein} de $\cal R$, not\'e $\pac{[\cal R]}$, est l'ensemble de tous les isomorphismes partiels de $X$ dont le graphe est contenu dans $\cal R$. Une application bor\'elienne $f : A \longrightarrow X$ est un {\it morphisme int\'erieur partiel} si tout \'el\'ement de~$A$ est $\cal R$-\'equivalent \`a son image par $f$. Si $\phi : A \rightarrow B$ est un \'el\'ement du pseudo-groupe plein de $\cal R$ et $\cal S$ une sous-relation de $\cal R$ d\'efinie sur $B$, on d\'efinit alors sur $A$ une sous-relation de $\cal R$ not\'ee $\phi^{-1}{\cal S}\phi$ : deux \'el\'ements $x$ et $y$ sont $\phi^{-1}{\cal S}\phi$-\'equivalents si par d\'efinition $\phi(x)$ et $\phi(y)$ sont $\cal S$-\'equivalents. Ainsi, $\cal S$ et $\phi^{-1}{\cal S}\phi$ sont des sous-relations isomorphes via $\phi$. On dit que $\phi^{-1}{\cal S}\phi$ est la sous-relation d\'eduite de $\cal S$ par {\it conjugaison} par $\phi$. Deux sous-relations $\cal S$ et ${\cal S}'$ de $\cal R$ d\'efinies sur les parties bor\'eliennes~$A$ et $A'$ de $X$ sont conjugu\'ees dans $\cal R$ si elles sont orbitalement \'equivalentes via un \'el\'ement $\phi : A \longrightarrow A'$ du pseudo-groupe plein de $\cal R$, autrement dit si ${\cal S}'=\phi^{-1}{\cal S}\phi$. Enfin, nous dirons que $\cal S$ et ${\cal S}'$ sont {\it stablement conjugu\'ees} dans $\cal R$ s'il existe des domaines complets $A$ et~$A'$ de $\cal S$ et ${\cal S}'$ respectivement sur lesquels les restrictions de $\cal S$ et ${\cal S}'$ sont conjugu\'ees dans $\cal R$.

\subsection{Espaces fibr\'es standards et actions}

Un {\it espace fibr\'e standard} $(F,{\cal B}_F,\pi)$ est la donn\'ee d'un espace bor\'elien standard $F$ sur $X$ et d'une application bor\'elienne (appel\'ee projection) $\pi: F \rightarrow X$ surjective \`a pr\'e-images d\'enombrables. La fibre $F_x$ d'un \'el\'ement $x$ de $X$ est la pr\'e-image de $x$ par $\pi$. Comme sous-ensemble bor\'elien de $X \times X$, la relation $\cal R$ d\'efinit naturellement deux espaces fibr\'es standards sur $X$ via les projections~$\pi_l$ et~$\pi_r$ respectivement sur la premi\`ere et deuxi\`eme coordonn\'ee. Une {\it section bor\'elienne}~$s$ de $F$ est une application bor\'elienne de $X$ dans $F$ telle que $\pi \circ s$ soit \'egale \`a l'identit\'e. Si $A$ est une partie bor\'elienne de $X$ et si $s$ n'est d\'efinie que sur $A$, alors nous parlerons de {\it section partielle}. Un espace fibr\'e standard sur $X$ admet toujours une section bor\'elienne. Ceci est une cons\'equence du th\'eor\`eme suivant (voir \cite{MR0217751}, \cite{MR1321597}) :

\begin{theoreme}[Th\'eor\`eme de s\'election]\label{selection}
Soit $F$ un espace fibr\'e standard sur $X$. Alors il existe une famille d\'enombrable de sections partielles de $X$ dont les images forment une partition (bor\'elienne et d\'enombrable) de $F$. De plus, on peut toujours supposer qu'au moins l'une de ces sections partielles est une section bor\'elienne, c'est-\`a-dire d\'efinie sur $X$ tout entier.
\end{theoreme}

\noindent
Voici trois applications imm\'ediates de ce th\'eor\`eme :

\begin{itemize}
\smallskip
\item si $F$ un espace fibr\'e standard sur $X$, alors il existe une num\'erotation bor\'elienne des fibres de $F$, c'est-\`a-dire une application bor\'elienne $N : F \longrightarrow \Nstar$ telle que la restriction de $N$ \`a toute fibre de $F$ soit injective. De plus, quitte \`a renum\'eroter les fibres de $F$, on peut toujours supposer que dans chaque fibre la num\'erotation commence \`a $1$ et ne saute pas d'entiers naturels ;
\item si $f$ est une r\'eduction de $\cal R$ \`a ${\cal R}'$, c'est-\`a-dire une application bor\'elienne $f : X \longrightarrow X'$ telle que deux \'el\'ements de $X$ sont $\cal R$-\'equivalents si et seulement si leurs images par $f$ sont ${\cal R}'$-\'equivalents (cf. \cite{MR1900547}), alors il existe un domaine complet $A$ de $\cal R$ tel que la restriction de~$f$ \`a~$A$ soit une \'equivalence orbitale de $\rest{\cal R}{A}$ \`a $\rest{{\cal R}'}{f(A)}$ ;
\item soit $\cal R$ une relation d'\'equivalence bor\'elienne sur $X$ et $A$ un domaine complet de $\cal R$. Alors il existe un morphisme int\'erieur de $\cal R$ d\'efinie sur $X$ et dont l'image est contenue dans $A$.
\end{itemize}

\bigskip

Nous allons maintenant introduire la notion d'{\it action} pour une relation d'\'equivalence bor\'elienne sur un espace fibr\'e standard $F$ sur $X$. Rappelons d'abord que le {\it produit fibr\'e} de deux espaces fibr\'es standards $(F',\pi')$ et $(F'',\pi'')$ sur $X$ est l'espace fibr\'e standard $(F,\pi)$ o\`u
\centers{$F = F' \star F''=\paa{(t',t'') \in F' \times F'' \ ; \ \pi'(t')=\pi''(t'')}$}
\noindent
et $\pi$ l'application bor\'elienne de $F$ dans $X$ d\'efinie par $\pi(t',t'')=\pi'(t')$.

\begin{definition}[Gaboriau, \cite{MR1953191}]\label{definition action}
Une {\it $\cal R$-action} (\`a gauche) sur l'espace fibr\'e standard $(F,\pi)$ sur~$X$ est une application bor\'elienne
\centers{$\application{({\cal R},\pi_r) \star (F,\pi)}{F}{\pa{(x,y),t}}{(x,y) \cdot t}$}
\noindent
telle que, pour tout triplet $(x,y,z)$ d'\'el\'ements $\cal R$-\'equivalents de $X$ et pour tout $t$ appartenant \`a $F$ dans la fibre de $z$, on ait
\centers{$(z,z) \cdot t=t \qqetqq (x,y) \cdot \pa{(y,z) \cdot t}=(x,z) \cdot t.$}
\noindent
On dit alors que $(F,\pi)$ est un {\it $\cal R$-espace fibr\'e standard} sur $X$ et que $\cal R$ agit sur $F$.
\end{definition}

\petiteremarque
Pour qu'elle ait un sens la formule du produit ci-dessus impose que $(x,y) \cdot t$ soit un \'el\'ement dans la fibre de $x$. De m\^eme on d\'efinit la notion d'action \`a droite que nous rencontrerons \'egalement par la suite.

\medskip

Soit $F$ un $\cal R$-espace fibr\'e standard sur $X$. L'{\it orbite} d'un \'el\'ement $f_x$ de $F$ dans la fibre de~$x$ de~$X$ est l'ensemble des $(y,x) \cdot f_x$ o\`u $y$ d\'ecrit la $\cal R$-classe de $x$. En particulier, l'action de $\cal R$ sur $F$ engendre une relation d'\'equivalence bor\'elienne not\'ee ${\cal R}_F$ sur $F$ : $f_x$ et $f_y$ sont ${\cal R}_F$-\'equivalents si, par d\'efinition, $(x,y) \cdot f_y=f_x$. Puisque deux \'el\'ements ${\cal R}_F$-\'equivalents de $F$ se projettent dans~$X$ sur des \'el\'ements $\cal R$-\'equivalents, la projection est donc un morphisme de relations d'\'equivalence bor\'eliennes.

\medskip

\exemplefondamental \label{espace fibre standard canonique}
$(F,\pi)=({\cal R},\pi_l)$ d\'efinit un ${\cal R}$-espace fibr\'e standard sur $X$ avec l'action \og horizontale \fg
\centers{$\application{({\cal R},\pi_r) \star ({\cal R},\pi_l)}{({\cal R},\pi_l)}{\pa{(x,y),(y,z)}}{(x,z).}$}

Les classes de ${\cal R}_F$ sont ici les fibres de $\pi_r : {\cal R} \longrightarrow X$. Nous dirons que $({\cal R},\pi_l)$ est le $\cal R$-{\it espace fibr\'e standard canonique} gauche associ\'e \`a $\cal R$.

\medskip

De la m\^eme fa\c{c}on, nous avons le $\cal R$-espace fibr\'e standard canonique droit $({\cal R},\pi_r)$ avec son action (\`a droite) \og verticale \fg

\centers{$\application{({\cal R},\pi_r) \star ({\cal R},\pi_l)}{({\cal R},\pi_r)}{\pa{(x,y),(y,z)}}{(x,z).}$}

\medskip

\petiteremarque
Si $(F,\pi)$ est un $\cal R$-espace fibr\'e standard sur $X$ et $A$ une partie bor\'elienne de $X$, on obtient alors une notion de $\rest{{\cal R}}{A}$-espace fibr\'e standard {\it induit} sur $A$ : il s'agit de la restriction de $(F,\pi)$ \`a $(\pi^{-1}(A),\rest{\pi}{\pi^{-1}(A)})$.

\medskip

Soit $s : A \longrightarrow F$ une section partielle d'un $\cal R$-espace fibr\'e standard $F$, le stabilisateur $\text{Stab}_{\cal R}(s)$ de $s$ est la sous-relation de $\cal R$ d\'efinie sur $A$ suivante : deux \'el\'ements $x$ et $y$ de $A$ sont $\text{Stab}_{\cal R}(s)$-\'equivalents si leurs images par $s$ sont ${\cal R}_F$-\'equivalents, autrement dit
\centers{$x \sim_{\text{Stab}_{\cal R}(s)}y \qquad \text{ssi} \qquad (y,x) \cdot s(x)=s(y).$}

Notons que la projection sur $X$ induit une \'equivalence orbitale entre la restriction de~${\cal R}_F$ \`a~$s(A)$ et $\text{Stab}_{\cal R}(s)$. Nous appellerons {\it satur\'e} de l'image de $s$ le ${\cal R}_F$-satur\'e de $s(A)$. Si $s$ est une section bor\'elienne de $F$ (c'est-\`a-dire si $A=X$), $s(X)$ est ${\cal R}_F$-satur\'e si et seulement si $\text{Stab}_{\cal R}(s)$ co\"incide avec~$\cal R$. Nous d\'efinissons maintenant la notion de morphisme entre $\cal R$-espaces fibr\'es standards sur~$X$.

\begin{definition}
Un morphisme $f : F \longrightarrow F'$ de $\cal R$-espaces fibr\'es standards est une application bor\'elienne de $F$ dans $F'$ telle que, si $x$ et $y$ sont $\cal R$-\'equivalents et si $f_x$ appartient \`a la fibre de $x$, alors
\centers{$f(F_x) \subset F'_x \quad \text{et} \quad  f\pa{(y,x) \cdot f_x}=(y,x) \cdot f(f_x).$}
\end{definition}

\petiteremarque
En particulier, un morphisme de $\cal R$-espaces fibr\'es standards de $F$ dans~$F'$ induit une application de la fibre $F_x$ dans la fibre $F'_x$ pour tout $x$ de $X$. Un tel morphisme est dit injectif (respectivement surjectif) si c'est une application bor\'elienne injective (resp. surjective) : au niveau des fibres, on obtient des applications de m\^eme nature. Remarquons \'egalement qu'un morphisme de $\cal R$-espaces fibr\'es standards de~$F$ dans $F'$ induit un morphisme de relations d'\'equivalence bor\'eliennes de~${\cal R}_F$ dans~${\cal R}_{F'}$. En particulier, un isomorphisme entre $F$ et $F'$ induit une \'equivalence orbitale entre~${\cal R}_F$ et ${\cal R}_{F'}$.

\bigskip

Soit $F$ un $\cal R$-espace fibr\'e standard sur $X$. L'action de $\cal R$ sur $F$ est {\it transitive} et $F$ est dit {\it homog\`ene} s'il existe un domaine complet de ${\cal R}_F$ sur lequel la restriction de la projection de $F$ sur $X$ est injective. Ainsi une action de $\cal R$ sur $F$ est transitive si et seulement s'il existe une section partielle~$s$ dite {\it saturante} d\'efinie sur une partie bor\'elienne $A$ de~$X$ dont l'image est un domaine complet de ${\cal R}_F$ : autrement dit, le satur\'e de $s(A)$ co\"{i}ncide avec $F$ ($A$ est n\'ecessairement un domaine complet de $\cal R$).

\begin{proposition}\label{stabilisateurs SOE}
Soit $s : A \longrightarrow F$ et $s' : A' \longrightarrow F$ deux sections partielles saturantes d'un $\cal R$-espace fibr\'e standard homog\`ene $F$. Alors leurs stabilisateurs $\text{Stab}_{\cal R}(s)$ et $\text{Stab}_{\cal R}(s')$ sont deux sous-relations de $\cal R$ stablement conjugu\'ees.
\end{proposition}

\demo
Consid\'erons le sous-ensemble bor\'elien $\Xi$ de $\cal R$ constitu\'e des couples d'\'el\'ements $(a,a')$ de $A \times A'$ tels que $s(a)$ et $s'(a')$ appartiennent \`a la m\^eme orbite. Puisque le satur\'e de $s'(A)$ contient $s(A)$, la restriction de $\pi_l$ \`a $\Xi$ induit sur $\Xi$ une structure d'espace fibr\'e standard. Soit $f'$ une section bor\'elienne de $(\Xi,\pi_l)$ : $\pi_r \circ f'$ est alors une r\'eduction $r$ de $\text{Stab}_{\cal R}(s)$ dans $A$ dans $\text{Stab}_{\cal R}(s')$ sur $A'$ qui de plus est un morphisme int\'erieur de $\cal R$ v\'erifiant
\centers{$s'(r(x))=(r(x),x) \cdot s(x).$}
Puisque le satur\'e de $s(A)$ contient le satur\'e de $s'(A')$, on en d\'eduit que $r(A)$ est un domaine complet de $\text{Stab}_{\cal R}(s')$ et ceci entra\^{i}ne que $r$ soit une \'equivalence orbitale stable entre $\text{Stab}_{\cal R}(s)$ et $\text{Stab}_{\cal R}(s')$.
\findemo

\medskip

\noindent
{\it Remarque 1 :}
La proposition pr\'ec\'edente assure, qu'\'etant donn\'e un $\cal R$-espace fibr\'e standard homog\`ene sur $X$, le stabilisateur d'une section partielle saturante est unique \`a \'equivalence orbitale stable pr\`es.

\smallskip

\noindent
{\it Remarque 2 :}
On d\'eduit de la proposition pr\'ec\'edente le fait suivant que nous utiliserons \`a plusieurs reprises : si $s''$ est une section partielle ({\it a priori} non saturante) de~$F$, alors le stabilisateur de $s'' : A'' \longrightarrow F$ est stablement conjugu\'e \`a une restriction de $\text{Stab}_{\cal R}(s)$ o\`u $s$ est une section partielle saturante. En effet, le satur\'e de $s''(A'')$ est un sous-espace fibr\'e standard $F''$ de la restriction de~$F$ \`a ${\cal R} \cdot A''$ qui est homog\`ene de section partielle saturante $s''$. Consid\'erons la projection sur $X$ de l'intersection de $F''$ et de $s(A)$ : c'est une partie bor\'elienne $B$ de $A$ telle que ${\cal R} \cdot B$ et ${\cal R} \cdot (A \setminus B)$ forment une partition de~$X$ et telle que la restriction de $s$ \`a $B$ soit alors une section partielle saturante de $F''$. La proposition pr\'ec\'edente assure alors que $\text{Stab}_{\cal R}(s'')$ et $\text{Stab}_{\cal R}(\rest{s}{B})=\rest{\text{Stab}_{\cal R}(s)}{B}$ soient stablement conjugu\'es.

\bigskip

Soit $F_1$ et $F_2$ deux $\cal R$-espaces fibr\'es standards sur $X$. D\'esignons par $s_1$ et $s_2$ des sections partielles de $F_1$ et $F_2$ d\'efinies sur la partie bor\'elienne $A$ de $X$. Supposons de plus que $f$ soit un morphisme de $\cal R$-espaces fibr\'es standards de $F_1$ dans $F_2$ qui envoie $s_1(A)$ sur $s_2(A)$ : le stabilisateur de $s_1$ est alors une sous-relation du stabilisateur de~$s_2$. Plus pr\'ecis\'ement, nous avons le lemme suivant :

\begin{lemme}\label{morphisme fibr\'es}
Soit $(F_1,\pi_1)$ et $(F_2,\pi_2)$ deux $\cal R$-espaces fibr\'es standards sur~$X$ et deux sections partielles $s_i : A \longrightarrow F_i$ d\'efinies sur un domaine complet $A$ de $\cal R$. Si~$s_1$ est saturante, alors il existe un morphisme $f : F_1 \longrightarrow F_2$ de $\cal R$-espaces fibr\'es standards qui envoie $s_1(A)$ sur $s_2(A)$ si et seulement si le stabilisateur $\text{Stab}_{\cal R}(s_1)$ est une sous-relation de $\text{Stab}_{\cal R}(s_2)$. Si de plus $s_2$ est saturante, alors~$f$ est surjectif.
\end{lemme}

\demo
Supposons que le stabilisateur $\text{Stab}_{\cal R}(s_1)$ soit une sous-relation de $\text{Stab}_{\cal R}(s_2)$. Si~$f_x$ appartient \`a la $F_1$-fibre d'un \'el\'ement $x$ de $X$, alors par hypoth\`ese il existe un \'el\'ement $y$ dans la classe de $x$ tel que $f_x$ appartienne \`a l'orbite de $s_1(y)$. On d\'efinit alors $f(f_x)$ comme l'image par le couple d'\'el\'ements $(x,y)$ de $s_2(y)$. Cette d\'efinition ne d\'epend pas du choix du repr\'esentant~$y$ puisque deux \'el\'ements $\text{Stab}_{\cal R}(s_1)$-\'equivalents sont $\text{Stab}_{\cal R}(s_2)$-\'equivalents par hypoth\`ese. Enfin, supposons que $s_2$ soit de plus saturante. On a alors
\centers{$F_2={\cal R} \cdot s_2(A)={\cal R} \cdot f(s_1(A))=f({\cal R} \cdot s_1(A))=f(F_1).$}
\findemo

\noindent
On en d\'eduit le fait important suivant :

\begin{proposition}\label{section domaine fondamental}
Soit $F$ un $\cal R$-espace fibr\'e standard sur $X$. S'il existe une section bor\'elienne $s$ de~$F$ telle que $s(X)$ soit un domaine fondamental de ${\cal R}_F$, alors $F$ est isomorphe au $\cal R$-espace fibr\'e standard canonique.
\end{proposition}

\demo
Le stabilisateur de $s$ \'etant la relation triviale par hypoth\`ese, le lemme pr\'ec\'edent assure l'existence d'un morphisme surjectif $f$ de $\cal R$-espaces fibr\'es standards entre $F$ et le $\cal R$-espace fibr\'e standard canonique (gauche) envoyant l'image de la section bor\'elienne $s$ sur la diagonale. Il ne reste plus qu'\`a voir que ce morphisme est injectif. Par l'absurde, supposons que $f_x$ et $f'_x$ soient deux \'el\'ements distincts dans la $F$-fibre d'un \'el\'ement $x$ de $X$ tels que leurs images par $f$ soient \'egales dans le $\cal R$-espace fibr\'e standard canonique. Il existerait alors deux \'el\'ements $y$ et $z$ distincts de $X$ tels que $f_x$ et $f'_x$ appartiennent respectivement aux orbites de $s(y)$ et $s(z)$. Mais par suite les images respectives par les couples d'\'el\'ements $(y,x)$ et $(z,x)$ des \'el\'ements $f(s(y))$ et $f(s(z))$ de la diagonale de $\cal R$ seraient \'egales dans le $\cal R$-espace fibr\'e standard canonique.
\findemo

\bigskip

Nous allons maintenant voir que certains $\cal R$-espaces fibr\'es standards peuvent se plonger dans le $\cal R$-espace fibr\'e standard canonique et nous utiliserons ce fait \`a plusieurs reprises.

\begin{lemme}\label{prolongement fibr\'e}
Soit $F$ un $\cal R$-espace fibr\'e standard sur $X$ admettant une section partielle $s$ d\'efinie sur$A$ telle que $s(A)$ soit un domaine fondamental de ${\cal R}_F$. Alors il existe un $\cal R$-espace fibr\'e standard $F'$ sur~$X$ contenant $F$ et ayant une section bor\'elienne $s'$ qui prolonge~$s$ \`a~$X$ et telle que $s'(X)$ soit un domaine fondamental pour l'action de $\cal R$ sur $F'$.
\end{lemme}

\petiteremarque
De l'existence d'une section bor\'elienne de $F'$ dont l'image est un domaine fondamental pour l'action de $\cal R$, on en d\'eduit que le $\cal R$-espace fibr\'e standard $F'$ est isomorphe au $\cal R$-espace fibr\'e standard canonique d'apr\`es la proposition pr\'ec\'edente.

\smallskip

\demo
Puisque l'image de $s$ est un domaine fondamental pour l'action de~$\cal R$ sur $F$, on en d\'eduit que $A$ est un domaine complet de $\cal R$. Consid\'erons l'espace fibr\'e standard r\'eunion disjointe de $\rest{F}{A}$ et de $\pi_l^{-1}(A) \cap \pi_r^{-1}(X \setminus A)$ qui est un espace fibr\'e standard~$F'$ sur $A$ via la restriction de~$\pi_l$. Par construction, $F'$ contient $\rest{F}{A}$ et les actions de $\rest{{\cal R}}{A}$ sur $\rest{F}{A}$ et de $\rest{{\cal R}}{A}$ sur la restriction \`a $A$ de l'espace fibr\'e standard canonique gauche se prolongent en une action de $\rest{{\cal R}}{A}$ sur $F'$. Enfin, puisque $A$ est un domaine complet de~$\cal R$, $F'$ s'\'etend par \'equivariance en un $\cal R$-espace fibr\'e standard sur~$X$ comme souhait\'e en d\'efinissant la section bor\'elienne $s'$ de $F'$ de telle sorte qu'elle co\"{i}ncide avec $s$ sur $A$ et avec \og l'identit\'e \fg sur le compl\'ementaire de~$A$ dans~$X$.
\findemo

\bigskip

Nous avons d\'ej\`a mentionn\'e que la donn\'ee d'un $\cal R$-espace fibr\'e standard $F$ d\'efinit naturellement une relation d'\'equivalence bor\'elienne ${\cal R}_F$ sur l'espace bor\'elien standard~$F$. Le cas o\`u ${\cal R}_F$ est lisse va particuli\`erement nous int\'eresser dans notre \'etude des relations d'\'equivalence bor\'eliennes arborables (cf. \S \ \ref{actions quasi-libres}).

\begin{definition}[Action lisse]
Soit $F$ un $\cal R$-espace fibr\'e standard sur $X$. L'action de $\cal R$ sur $F$ est lisse (on dit aussi que $\cal R$ agit de mani\`ere lisse sur $F$) si la relation d'\'equivalence bor\'elienne ${\cal R}_F$ sur~$F$ est lisse.
\end{definition}

\exemple
L'action de $\cal R$ sur l'espace fibr\'e standard canonique est lisse ; il en est de m\^eme de toute sous-relation de $\cal R$ puisqu'une sous-relation d'une relation lisse est elle-m\^eme lisse. Plus g\'en\'eralement, si $\cal R$ agit de mani\`ere lisse sur un espace fibr\'e standard $F$, il en est de m\^eme de chacune de ses sous-relations.

\medskip

Nous allons donner une caract\'erisation des actions lisses que nous utiliserons constamment et qui est une cons\'equence du lemme suivant.

\begin{lemme}\label{partition}
Pour tout $\cal R$-espace fibr\'e standard $F$ sur $X$, il existe une famille d\'enombrable $(s_i : A_i \longrightarrow F)_{i \in I}$ de sections partielles de $F$ dont les satur\'es $F_i$ des images $s_i(A_i)$ forment une partition ${\cal R}_F$-invariante de $F$.
\end{lemme}

\petiteremarque
Puisque par d\'efinition $\rest{{\cal R}_F}{F_i}$ et $\rest{{\cal R}_F}{s_i(A_i)}$ sont stablement orbitalement \'equivalentes, on en d\'eduit qu'il en est de m\^eme de $\rest{{\cal R}_F}{F_i}$ et $\text{Stab}_{\cal R}(s_i)$.

\smallskip

\demo
Donnons-nous une num\'erotation bor\'elienne des fibres de $F$ et consid\'erons la section bor\'elienne qui, \`a chaque \'el\'ement de $X$, associe le plus petit \'el\'ement dans sa fibre. Si le compl\'ementaire du satur\'e de l'image de cette section bor\'elienne est vide, c'est termin\'e. Sinon on consid\`ere la section partielle d\'efinie par les plus petits \'el\'ements restants dans chaque fibre puis le satur\'e de l'image de cette derni\`ere. En continuant ainsi, on construit \`a chaque \'etape une nouvelle section partielle en prenant dans chaque fibre les \'el\'ements les plus petits restants dans le compl\'ementaire des satur\'es des images des sections partielles pr\'ec\'edemment construites. Cette construction fournit une exhaustion de $F$ puisqu'un \'el\'ement de num\'ero $n$ dans une fibre de $F$ a forc\'ement \'et\'e consid\'er\'e avant la $n\ieme$ \'etape.
\findemo

\begin{definition}[Action quasi-libre]\label{definition quasi-libre}
Soit $F$ un $\cal R$-espace fibr\'e standard sur $X$. L'action de $\cal R$ sur $F$ est quasi-libre (on dit encore que $\cal R$ agit quasi-librement sur $F$) si le stabilisateur de toute section partielle de $F$ est une sous-relation lisse de $\cal R$.
\end{definition}

On obtient alors la caract\'erisation suivante :

\begin{proposition}\label{quasi-libre=lisse}
\'Etant donn\'e un $\cal R$-espace fibr\'e standard $F$ sur $X$, l'action est quasi-libre si et seulement si $\cal R$ agit de mani\`ere lisse sur $F$, c'est-\`a-dire si ${\cal R}_F$ est lisse.
\end{proposition}

\demo
L'implication r\'eciproque est claire car nous avons d\'ej\`a mentionn\'e que le stabilisateur de toute section partielle $s : A \longrightarrow F$ est orbitalement \'equivalent \`a $\rest{{\cal R}_F}{s(A)}$ qui est une sous-relation de ${\cal R}_F$. Supposons maintenant l'action quasi-libre et construisons un domaine fondamental de ${\cal R}_F$. Le lemme pr\'ec\'edent assure l'existence d'une famille d\'enombrable $(s_i)_{i \in I}$ de sections partielles de $F$ dont les satur\'es des images forment une partition ${\cal R}_F$-invariante de $F$. Le stabilisateur de chacune de ces sections partielles $s_i$ \'etant lisse, consid\'erons la restriction de $s_i$ \`a un domaine fondamental $D_i$ de son stabilisateur : la r\'eunion sur $I$ des $s(D_i)$ est un domaine fondamental de~${\cal R}_F$.
\findemo

\bigskip

\label{exemple fondamental fibre relations}
Nous allons \`a pr\'esent introduire une classe de $\cal R$-espaces fibr\'es standards fondamentaux pour une relation d'\'equivalence bor\'elienne $\cal R$ sur $X$. Si $\cal S$ est une sous-relation de $\cal R$ d\'efinie sur $X$, alors $\cal S$ agit sur $\cal R$ via l'action induite par celle de \og $\cal R$ sur $\cal R$ \fg. En effet, consid\'erons $\cal R$ munie de ses deux structures de $\cal R$-espace fibr\'e standard canoniques et remarquons que la projection $\pi_l : {\cal R} \longrightarrow X$ est invariante sous l'action verticale de $\cal R$ (et donc de $\cal S$) sur $({\cal R},\pi_r)$. De plus, puisque $\cal S$ est une sous-relation de $\cal R$, elle agit \'egalement de mani\`ere lisse sur $({\cal R},\pi_r)$. Soit $\tf{\cal R}{\cal S}$ l'espace quotient de $\cal R$ par la relation d'\'equivalence bor\'elienne engendr\'ee par cette action. Comme les actions horizontale et verticale de \og $\cal R$ sur~$\cal R$ \fg commutent, on en d\'eduit que c'est un $\cal R$-espace fibr\'e standard sur~$X$ dont la projection sur $X$ et l'action de $\cal R$ sont induites par celles du $\cal R$-espace fibr\'e standard canonique gauche $({\cal R},\pi_l)$ : c'est le {\it $\cal R$-espace fibr\'e standard canonique gauche associ\'e au couple $({\cal R},{\cal S})$}. De plus, la diagonale $d$ de ce dernier passe au quotient sous l'action de $\cal S$ et le $\cal R$-espace fibr\'e standard $(\tf{{\cal R}}{{\cal S}},\pi_l)$ est naturellement muni d'une section bor\'elienne $d_{\cal S}$ dont l'image est un domaine complet de ${\cal R}_{\tf{\cal R}{\cal S}}$ et dont le stabilisateur est $\cal S$. Notons \'egalement que dans le cas o\`u la sous-relation $\cal S$ est triviale, $(\tf{{\cal R}}{{\cal S}},\pi_l)$ n'est autre que le $\cal R$-espace fibr\'e standard canonique gauche.

\smallskip

\petiteremarque
La construction pr\'ec\'edente s'\'etend au cas de sous-relations d\'efinies sur un domaine complet $A$ de $\cal R$ en consid\'erant l'action de $\cal S$ sur $({\cal R} \cap \pi_r^{-1}(A),\pi_r)$. Dans ce cas, le $\cal R$-espace fibr\'e standard $\tf{\cal R}{\cal S}$ est naturellement muni d'une section partielle $d_{\cal S}$ d\'efinie sur $A$ dont l'image est un domaine complet de  ${\cal R}_{\tf{\cal R}{\cal S}}$ et dont le stabilisateur est $\cal S$. En particulier, $\tf{\cal R}{\cal S}$ est un $\cal R$-espace fibr\'e standard homog\`ene et~$d_{\cal S}$ une section partielle saturante.

\medskip

Remarquons que l'on peut donner une description explicite de $\tf{\cal R}{\cal S}$ en tant qu'espace fibr\'e standard sur $X$. En effet, il suffit pour ceci de se donner une num\'erotation bor\'elienne des fibres de l'espace fibr\'e standard canonique gauche. Consid\'erons la partie bor\'elienne de $({\cal R},\pi_l)$ constitu\'ee des paires $(x,y)$ o\`u $y$ d\'esigne l'\'el\'ement de plus petit num\'ero dans sa $\cal S$-classe. L'espace fibr\'e standard obtenu est alors isomorphe \`a l'espace fibr\'e standard $\tf{\cal R}{\cal S}$.

\medskip
\label{fibre a gauche=fibre a droite pour relations}
De la m\^eme mani\`ere, on peut \'egalement consid\'erer le $\cal R$-espace fibr\'e standard canonique droit et faire agir $\cal S$ \`a gauche. On obtient alors le $\cal R$-espace fibr\'e standard~$\tfi{{\cal S}}{{\cal R}}$. Soit $S$ la sym\'etrie par rapport \`a la diagonale :
\centers{$\fonction{S}{{\cal R}}{{\cal R}}{(x,y)}{(y,x).}$}
Soit $x$, $x'$ et $y$ trois \'el\'ements $\cal R$-\'equivalents tels que $x$ et $x'$ appartiennent \`a la m\^eme $\cal S$-classe. Puisque les images de $S(x,y)$ et $S(x',y)$ sont \'egales dans l'espace quotient $\tfi{{\cal S}}{{\cal R}}$, cette application passe au quotient et on obtient l'application bor\'elienne
\centers{$\tilde{S} : \pa{\tf{{\cal R}}{{\cal S}},\pi_l} \longrightarrow \pa{\tfi{{\cal S}}{{\cal R}},\pi_r}$}
\noindent
qui est un isomorphisme de $\cal R$-espaces fibr\'es standards.

\subsection{Actions quasi-libres et arboralit\'e}\label{actions quasi-libres}

\label{definition arborage}
Rappelons qu'un (L-)graphage (cf. \cite{MR1366313}, \cite{MR2095154}) de $\cal R$ munit canoniquement chaque classe de $\cal R$ d'une structure de graphe connexe, son graphe de Cayley, dont les sommets sont les \'el\'ements de cette classe. Un (L-)graphage de $\cal R$ est un {\it (L-)arborage} si les graphes de Cayley de chaque classe de $\cal R$ sont des arbres. Une relation d'\'equivalence bor\'elienne est dite {\it arborable} si elle admet un (L-)arborage. Les relations d'\'equivalence bor\'eliennes arborables jouent un r\^ole central dans la classe des relations d'\'equivalence bor\'eliennes, au m\^eme titre que les groupes libres dans la classe des groupes. Bien souvent, les invariants qui ont \'et\'e introduits pour \'etudier les relations d'\'equivalence mesur\'ees (voir par exemple \cite{MR1646912}, \cite{MR1953191}) ont d'abord \'et\'e calcul\'es pour les relations arborables. Dans ce paragraphe, nous allons donner une caract\'erisation des relations d'\'equivalence bor\'eliennes arborables (th. \ref{quasi-libre=arborable}) analogue \`a celle concernant les groupes libres : {\it un groupe est libre si et seulement s'il agit librement sur un arbre} (voir par exemple \cite{MR0476875}).

\medskip

Comme nous l'avons d\'ej\`a vu, une relation d'\'equivalence bor\'elienne $\cal R$ sur $X$ d\'efinit canoniquement un $\cal R$-espace fibr\'e standard sur $X$ via la projection $\pi_l$. Supposons que $\cal R$ soit arborable et d\'esignons par $\Phi$ un (L-)arborage de $\cal R$. Dans ce cas, les fibres de l'espace fibr\'e standard canonique gauche sont naturellement munies d'une structure d'arbre (cf. {\it infra}). D\'efinissons \`a pr\'esent les $\cal R$-champs de graphes bor\'eliens et nous nous int\'eressons plus particuli\`erement au cas des $\cal R$-champs d'arbres bor\'eliens qui sont les acteurs principaux de ce travail.

\begin{definition}[$\cal R$-arboretum]\label{definition arboretum}
Un $\cal R$-{\it champ de graphes bor\'elien} $({\cal A},\pi)$ sur $X$ est un graphe dont les espaces de sommets et d'ar\^etes sont des $\cal R$-espaces fibr\'es standards $({\cal A}^0,{{\pi}^0})$ et $({\cal A}^1,{{\pi}^1})$ sur $X$ et tel que les applications sommet origine $o : {\cal A}^1 \longrightarrow {\cal A}^0$, sommet terminal $t : {\cal A}^1 \longrightarrow {\cal A}^0$ et ar\^ete oppos\'ee $\bar{} : {\cal A}^1 \longrightarrow {\cal A}^1$ soient des morphismes de $\cal R$-espaces fibr\'es standards.

La fibre dans $\cal A$ d'un \'el\'ement de $x$, not\'ee ${\cal A}_x$, est le sous-graphe d'ensemble de sommets $({\pi}^0)^{-1}(x)$ et d'ensemble d'ar\^etes $({\pi}^1)^{-1}(x)$. Si ${\cal A}_x$ est un arbre pour tout $x$ de $X$, nous dirons que $({\cal A},\pi)$ est un $\cal R$-arboretum.
\end{definition}

\petiteremarque
Les d\'efinitions pr\'ec\'edentes dans le cas de la relation d'\'equivalence bor\'elienne \`a classes triviales sur $X$ permettent de d\'efinir les notions de {\it champ de graphes bor\'elien} et d'{\it arboretum} sur~$X$.

\smallskip

\exemplefondamental
Tout graphage $\Phi$ (\cite{MR2095154}) sur $X$ d\'efinit un ${\cal R}_{\Phi}$-champ de graphes bor\'elien sur $X$ o\`u ${\cal R}_{\Phi}$ est la relation d'\'equivalence bor\'elienne engendr\'ee par~$\Phi$. L'espace des sommets est le ${\cal R}_{\Phi}$-espace fibr\'e standard canonique gauche et~${\cal R}_{\Phi}$ agit naturellement sur l'espace des ar\^etes d\'efini par le graphage $\Phi$ : l'image par le couple d'\'el\'ements ${\cal R}_{\Phi}$-\'equivalents $(x',x)$ de l'ar\^ete $(y,z)_x$ dans la fibre de $x$ est l'ar\^ete $(y,z)_{x'}$ dans la fibre de $x'$. Notons que l'action est, par d\'efinition, lisse sur l'espace des sommets et que la projection $\pi_r$ de l'espace des sommets dans $X$ est invariante sous l'action de ${\cal R}_{\Phi}$. Si $\Phi$ est un arborage, nous d\'esignerons par $({\cal A}_{\Phi},\pi)$ le ${\cal R}_{\Phi}$-arboretum canonique associ\'e \`a $\Phi$ sur $X$ ci-dessus.

\bigskip

Nous venons de voir dans l'exemple pr\'ec\'edent que si $\cal R$ est arborable, alors elle agit quasi-librement sur l'espace des sommets d'un arboretum. Puisque le stabilisateur d'une section partielle d'ar\^etes est toujours une sous-relation du stabilisateur de la section de sommets origines associ\'ee, on en d\'eduit que si $\cal R$ agit quasi-librement sur l'espace des sommets d'un arboretum (et plus g\'en\'eralement d'un champ de graphes bor\'elien), alors $\cal R$ agit \'egalement quasi-librement sur l'espace des ar\^etes. Ainsi, nous dirons que $\cal R$ agit {\it quasi-librement} sur un champ de graphes bor\'elien si~$\cal R$ agit quasi-librement sur l'espace des sommets.

\medskip

Nous allons \`a pr\'esent voir qu'il s'agit en fait d'une caract\'erisation des relations d'\'equivalence bor\'eliennes arborables et en donner un corollaire imm\'ediat (le reste de cette section est consacr\'e \`a la d\'emonstration du th\'eor\`eme).

\begin{theoreme}\label{quasi-libre=arborable}
Une relation d'\'equivalence bor\'elienne $\cal R$ est arborable si et seulement s'il existe une action quasi-libre de $\cal R$ sur un arboretum.
\end{theoreme}

Puisqu'une sous-relation d'une relation d'\'equivalence bor\'elienne lisse est lisse, la proposition~\ref{quasi-libre=lisse} et le th\'eor\`eme pr\'ec\'edent nous permettent de d\'emontrer un analogue du th\'eor\`eme de Nielsen-Schreier en th\'eorie des groupes : {\it un sous-groupe d'un groupe libre est libre}.

\begin{corollaire}[voir aussi \cite{MR1900547} et \cite{MR1728876}]\label{sous-relation d'une arborable}
Une sous-relation d'une relation d'\'equivalence bor\'elienne arborable est arborable.
\end{corollaire}

\petiteremarque
Si $A$ est une partie bor\'elienne de $X$ et si $\cal R$ est arborable, on d\'eduit du corollaire pr\'ec\'edent que la restriction de $\cal R$ \`a $A$ est encore arborable (il suffit de prolonger $\rest{\cal R}{A}$ par la relation triviale en dehors de $A$).

\medskip

\'Etant donn\'e une action quasi-libre d'une relation d'\'equivalence bor\'elienne $\cal R$ sur un arboretum~$\cal A$, nous allons d'abord montrer l'existence d'un sous-arboretum~${\cal A}'$ de $\cal A$ d\'efini sur une partie bor\'elienne de $X$ dont l'espace des sommets ${\cal A}'^0$ est un domaine fondamental de ${\cal R}_{{\cal A}^0}$. Nous nous ramenons ensuite au cas d'une action de~$\cal R$ sur un arboretum ${\cal A}^{''}$ ayant une section partielle de sommets $s : A \subset X \longrightarrow {\cal A}''^0$ dont l'image $s(A)$ est un domaine fondamental de ${\cal R}_{{\cal A}''^0}$, ce qui nous permettra de conclure gr\^ace au lemme suivant :

\begin{lemme}\label{arborable}
Si $\cal A$ est un $\cal R$-arboretum sur $X$ et $s$ une section partielle de sommets dont l'image est un domaine fondamental de ${\cal R}_{{\cal A}^0}$, alors $\cal R$ est arborable.
\end{lemme}

\demo
Commen\c{c}ons par supposer que $s$ soit une section bor\'elienne, c'est-\`a-dire d\'efinie sur $X$ tout entier. Comme nous l'avons d\'ej\`a vu (cf. prop. \ref{section domaine fondamental}), l'espace des sommets ${\cal A}^0$ s'identifie naturellement avec le $\cal R$-espace fibr\'e standard canonique gauche. L'ensemble des ar\^etes de $\cal A$ dont l'un des sommets appartient \`a l'image de la section bor\'elienne de sommets $s$ s'identifie alors \`a une partie bor\'elienne $Gr$ de~$X \times X$ qui, par construction, est un arborage de $\cal R$ : un couple d'\'el\'ements distincts $\cal R$-\'equivalents $(x,y)$ appartient \`a $Gr$ si, par d\'efinition, les sommets $(x,x)$ et $(x,y)$ sont adjacents dans l'arboretum $\cal A$.

\smallskip

Le cas g\'en\'eral s'en d\'eduit facilement gr\^ace au fait g\'en\'eral suivant. \'Etant donn\'e une relation d'\'equivalence bor\'elienne ${\cal R}'$ sur $X$, un domaine complet $A'$ de ${\cal R}'$ et un morphisme int\'erieur $r : X \longrightarrow A'$ qui co\"{i}ncide avec l'identit\'e sur $A'$, d\'efinissons la sous-relation ${\cal T}'$ de ${\cal R}'$ dont les classes sont les $r^{-1}(r(x))$ pour tout $x$ de~$X$ : en particulier, ${\cal T}'$ est lisse de domaine fondamental $A'$, donc arborable. Or, par construction, ${\cal R}'$ est le produit libre (cf. d\'ef. \ref{d\'efinition produit libre}) des sous-relations ${\cal T}'$ et ${\cal R}''$ o\`u~${\cal R}''$ co\"{i}ncide avec $\rest{{\cal R}'}{A'}$ sur $A'$ et avec la relation triviale sur le compl\'ementaire de $A'$. Enfin, le produit libre de deux relations d'\'equivalence bor\'eliennes arborables \'etant arborable (cf. prop. \ref{produit libre arborable}, voir aussi \cite{MR1900547} ou \cite{MR1728876}), on en d\'eduit que si $\rest{{\cal R}'}{A'}$ est arborable, alors il en est de m\^eme de~${\cal R}'$.
\findemo

\bigskip

Nous en venons au point central dans la d\'emonstration du th\'eor\`eme \ref{quasi-libre=arborable}.

\begin{proposition}\label{arbre d'isomorphismes partiels}
Si $\cal R$ agit quasi-librement sur un arboretum $\cal A$, alors il existe un sous-arboretum~${\cal A}'$ de la restriction de $\cal A$ \`a une partie bor\'elienne $A$ de $X$ tel que l'espace des sommets ${\cal A}'^0$ soit un domaine fondamental de ${\cal R}_{{\cal A}^0}$.
\end{proposition}

\petiteremarque
Notons que la proposition pr\'ec\'edente s'\'etend sans difficult\'e suppl\'ementaire au cas d'un $\cal R$-champ de graphes bor\'elien $\cal A$.

\smallskip

\demo
Donnons-nous une num\'erotation bor\'elienne index\'ee par les entiers dans chaque fibre de l'espace des sommets ${\cal A}^0$. Nous allons construire la partie bor\'elienne ${\cal A}'^0$ par \'etape \`a partir de sections partielles de sommets de ${\cal A}^0$. Soit $s_1$ la section bor\'elienne de sommets de ${\cal A}^0$ correspondant au num\'ero le plus petit dans chaque fibre. Puisque l'action de $\cal R$ sur $\cal A$ est quasi-libre, le stabilisateur de cette section est lisse. D\'esignons par $X_1$ un domaine fondamental de $\text{Stab}_{\cal R}(s_1)$ : cette partie bor\'elienne de $X$ est un domaine complet de $\cal R$ et par suite, le ${\cal R}_{{\cal A}^0}$-satur\'e de $U_1=s_1(X_1)$ rencontre toutes les fibres de ${\cal A}^0$. Notons $C_1$ le compl\'ementaire dans ${\cal A}^0$ du satur\'e de $U_1$. Si $C_1$ est vide, alors $U_1$ est un domaine fondamental pour l'action de $\cal R$ sur ${\cal A}^0$ et le r\'esultat est d\'emontr\'e avec ${\cal A}^{'0}=U_1$ et ${\cal A}^{'1}=\emptyset$.

\smallskip

Soit un entier naturel $n \geq 2$ . Supposons avoir construit les familles finies de sections partielles $S_j=\paa{s_i \ ; \ i \in I_j}$ ($1 \leq j \leq n-1$) de ${\cal A}^0$ dont les domaines de d\'efinition sont des parties bor\'eliennes de $X_1$ et v\'erifiant les deux conditions suivantes : la r\'eunion $U_{n-1}$ des images de ces sections partielles est une partie bor\'elienne de ${\cal A}^0$ sur laquelle la restriction de ${\cal R}_{{\cal A}^0}$ est triviale et l'intersection de $U_{n-1}$ avec chaque fibre de $\cal A$ au-dessus de $X_1$ d\'efinit une partie connexe. Notons $C_{n-1}$ le compl\'ementaire dans ${\cal A}^0$ du satur\'e de $U_{n-1}$. Si ce dernier est vide, la proposition est d\'emontr\'ee avec ${\cal A}^{'0}=U_{n-1}$ et ${\cal A}^{'1}$ la partie bor\'elienne de~${\cal A}^1$ constitu\'ee des ar\^etes dont les deux extr\'emit\'es appartiennent \`a~$U_{n-1}$.

Sinon, il existe des \'el\'ements de $C_{n-1}$ \`a distance unit\'e de $U_{n-1}$. En effet, les fibres de $\cal A$ \'etant connexes, il existe des \'el\'ements de $C_{n-1}$ \`a distance unit\'e du satur\'e de $U_{n-1}$ et $\cal R$ agit en pr\'eservant les distances sur ${\cal A}^0$. D\'esignons par $s_n$ la section partielle (d\'efinie sur une partie bor\'elienne de $X_1$ a priori) de sommets de ${\cal A}^0$ qui a pour image le plus petit \'el\'ement \`a distance unit\'e de $U_{n-1}$ dans chaque fibre. L'action \'etant quasi-libre, notons $X_n$ un domaine fondamental du stabilisateur de $s_n$ : c'est une partie bor\'elienne de $X_1$. Quitte \`a d\'ecouper $X_n$ en un nombre fini de parties bor\'eliennes, on peut raffiner $s_n$ en un nombre fini $S_n=\paa{s_i \ ; \ i \in I_n}$ de sections partielles dont chacune est telle que les \'el\'ements de son image soient \`a distance unit\'e de l'image de l'une des sections partielles d\'ej\`a construites.

\smallskip

Montrons d\'esormais que le satur\'e de la r\'eunion des $U_k$ est \'egal \`a l'espace des sommets ${\cal A}^0$. Sinon il existerait un \'el\'ement $a$ appartenant au compl\'ementaire du satur\'e de la r\'eunion des $U_k$ dans ${\cal A}^0$ dont la projection $\pi(a)$ serait un \'el\'ement de $X_1$. On peut supposer que dans la fibre de $\pi(a)$, le sommet $a$ soit \`a distance unit\'e de la partie connexe form\'ee par la r\'eunion des sommets $s_k(\pi(a))$ et de num\'ero minimal. D\'esignons par $n_0$ l'entier naturel tel que les sommets $a$ et $s_{n_0}(\pi(a))$ soient adjacents. Mais alors, par construction de la section partielle $s_{n_0+1}$, le sommet~$a$ devrait \^etre l'image de $\pi(a)$ par $s_{n_0+1}$, ce qui est absurde.

\smallskip

Il ne reste plus qu'\`a poser ${\cal A}^{'0}=\union{k \geq 1}{} U_k$ et ${\cal A}^{'1}$ la partie bor\'elienne de ${\cal A}^1$ constitu\'ee des ar\^etes dont les deux extr\'emit\'es appartiennent \`a ${\cal A}^{'0}$. Par construction, deux sommets de $\cal A$ dans la m\^eme orbite n'appartiennent pas tous les deux \`a ${\cal A}^{'0}$ ; puisque le satur\'e de la r\'eunion des $U_k$ est \'egal \`a l'espace des sommets ${\cal A}^0$, on en d\'eduit que~${\cal A}^{'0}$ est un domaine fondamental de ${\cal R}_{{\cal A}^0}$. Enfin, notre construction assure \`a chaque \'etape que les fibres de ${\cal A}'$ soient connexes.
\findemo

\medskip

\`A partir du sous-arboretum ${\cal A}^{'}$ construit dans la proposition pr\'ec\'edente, nous allons effectuer une op\'eration de contraction dans les fibres et obtenir un nouvel $\cal R$-arboretum ${\cal A}^{''}$ ayant une section partielle de sommets dont l'image est un domaine fondamental de ${\cal R}_{{\cal A}''^0}$.

\begin{lemme}\label{contraction}
Il existe un $\cal R$-arboretum ${\cal A}^{''}$ sur $X$ ayant une section partielle de sommets $s$ d\'efinie sur $A$ dont l'image est un domaine fondamental de ${\cal R}_{{\cal A}''^0}$.
\end{lemme}

\demo
Soit $s$ une section bor\'elienne de sommets de ${\cal A}'$ : son image est un domaine fondamental de la relation d'\'equivalence bor\'elienne lisse $\cal S$ d\'efinie sur ${\cal A}'^0$ dont les classes sont les fibres de ${\cal A}'^0$. Pour tout \'el\'ement $x$ de $A$, contractons le sous-arbre ${\cal A}_x'$ sur le sommet $s(x)$ dans l'arbre ${\cal A}_x$ (l'espace quotient est un espace bor\'elien standard isomorphe \`a $s(A)$ par construction). On prolonge cette op\'eration par $\cal R$-\'equivariance et on obtient ainsi un $\cal R$-arboretum ${\cal A}^{''}$ sur $X$ ayant une section partielle de sommets privil\'egi\'ee (encore not\'ee $s$) d\'efinie sur $A$. Si $x$ et $y$ sont deux \'el\'ements $\cal R$-\'equivalents de $A$, l'image par le couple d'\'el\'ements $(x,y)$ du sous-arbre ${\cal A}_y'$ de la fibre de $y$ dans la fibre de $x$ est un sous-arbre disjoint de ${\cal A}_x'$ puisque par hypoth\`ese ${\cal A}'^0$ est un domaine fondamental de l'action de $\cal R$ sur l'espace des sommets ${\cal A}^0$ de $\cal A$. Par suite, $s$ est une section partielle dont l'image est un domaine fondamental de ${\cal R}_{{\cal A}''^0}$.
\findemo

\medskip

Ceci termine la d\'emonstration du th\'eor\`eme \ref{quasi-libre=arborable} dont nous rappelons bri\`evement les trois \'etapes de la r\'eciproque. \'Etant donn\'e une action quasi-libre de $\cal R$ sur un arboretum $\cal A$ sur $X$ :
\begin{enumerate}
\item construction d'un sous-arboretum ${\cal A}'$ d\'efini sur une partie bor\'elienne $A \subset X$ de $\rest{\cal A}{A}$ dont l'espace des sommets ${\cal A}'^0$ est un domaine fondamental de ${\cal R}_{{\cal A}^0}$ (prop. \ref{arbre d'isomorphismes partiels}) ;
\item op\'eration de contraction pour se ramener au cas d'une action de $\cal R$ sur un arboretum ${\cal A}^{''}$ ayant une section partielle de sommets $s$ dont l'image $s(A)$ est un domaine fondamental de ${\cal R}_{{\cal A}''^0}$ (lem. \ref{contraction}) ;
\item une telle relation est alors arborable (lem. \ref{arborable}).
\end{enumerate}

\section{$\cal R$-arboretums et amalgames de sous-relations}
\label{sec:amalgames de sous-relations}

Nous allons nous int\'eresser dans cette partie aux relations d'\'equivalence bor\'eliennes qui sont des produits libres de deux sous-relations et, plus g\'en\'eralement, des produits amalgam\'es suivant une sous-relation commune. Il s'agit en fait de discuter des liens \'etroits entre d\'ecomposition d'une telle relation d'\'equivalence bor\'elienne et ses actions sur des arboretums. Nous avons vu (cf. \S \ \ref{actions quasi-libres}) que si une relation d'\'equivalence bor\'elienne agit quasi-librement sur un arboretum, alors elle est arborable. Nous souhaitons nous affranchir progressivement de l'hypoth\`ese faite sur l'action et nous traiterons le cas g\'en\'eral dans le prochain paragraphe (cf. \S \ \ref{cas g\'en\'eral}). Cette partie a pour but de d\'ecrire la situation dans le cas d'un \og espace quotient \fg simple.

\subsection{Produit libre de deux sous-relations}

Soit $\cal R$ une relation d'\'equivalence bor\'elienne engendr\'ee par deux sous-relations~${\cal R}_1$ et ${\cal R}_2$ sur~$X$. Autrement dit, $\cal R$ est la plus petite relation d'\'equivalence dont les classes contiennent celles de ${\cal R}_1$ et ${\cal R}_2$. Soit un entier naturel $n \geq 2$. Un $n$-uplet $\liste{x}$ d'\'el\'ements de $X$ est dit {\it r\'eduit} si :

\begin{itemize}
\item pour tout $i$ dans $\intn{1}{n-1}$, le couple $(x_i, x_{i+1})$ est un couple d'\'el\'ements ${\cal R}_1$ ou ${\cal R}_2$-\'equivalents et deux couples successifs ne sont pas ${\cal R}_j$-\'equivalents pour le m\^eme $j$ ;
\item $x_1 \neq x_2 \qsiq n=2$.
\end{itemize}

\begin{definition}[Gaboriau, \cite{MR1728876}]\label{d\'efinition produit libre}
La relation d'\'equivalence bor\'elienne $\cal R$ sur~$X$ est le produit libre des sous-relations ${\cal R}_1$ et ${\cal R}_2$ si $\cal R$ est engendr\'ee par ${\cal R}_1$ et ${\cal R}_2$ et si, pour tout $n$-uplet r\'eduit $\liste{x}$ d'\'el\'ements de $X$, $x_n$ est diff\'erent de $x_1$. Dans ce cas, on note ${\cal R}={\cal R}_1 \star {\cal R}_2$.
\end{definition}

\exemple
Soit $\Gamma_1$ et $\Gamma_2$ deux groupes d\'enombrables et une action libre $\alpha$ du produit libre $\Gamma=\Gamma_1 \star \Gamma_2$ sur $X$. La relation d'\'equivalence bor\'elienne ${\cal R}_{\alpha}$ dont les classes sont les orbites de l'action de $\Gamma$ est le produit libre des sous-relations ${\cal R}_{\alpha_1}$ et ${\cal R}_{\alpha_2}$ dont les classes sont les orbites des restrictions de $\alpha$ \`a $\Gamma_1$ et $\Gamma_2$.

\medskip

Avec les techniques que nous avons d\'evelopp\'ees jusqu'\`a pr\'esent, il est tr\`es facile de voir que l'arborabilit\'e est stable par produit libre.

\begin{proposition}\label{produit libre arborable}
Soit $\cal R$ une relation d'\'equivalence bor\'elienne sur $X$. Si $\cal R$ est le produit libre de deux sous-relations arborables, alors $\cal R$ est arborable.
\end{proposition}

\demo
Il suffit de construire un $\cal R$-arboretum tel que l'espace des sommets soit le $\cal R$-espace fibr\'e standard canonique (cas particulier d'une section bor\'elienne d\'efinie sur $X$ dans le lemme \ref{arborable}). \'Etant donn\'e des arborages de~${\cal R}_1$ et~${\cal R}_2$ respectivement, on s'int\'eresse aux arboretums canoniques associ\'es \`a ces arborages (cf. ex. fond. suivant d\'ef. \ref{definition arboretum}). Rappelons que les espaces des sommets de ces arboretums sont respectivement les ${\cal R}_1$ et ${\cal R}_2$-espaces fibr\'es standards canoniques gauches. Consid\'erons la r\'eunion disjointe de ces deux arboretums sur $X$, quotient\'ee par la relation d'\'equivalence qui identifie $d_1(x)$ et $d_2(x)$ pour tout \'el\'ement $x$ de $X$. On obtient ainsi un nouvel arboretum ${\cal A}'$ sur $X$ dont on note $d$ la section bor\'elienne commune $d_1=d_2$. Puisque ${\cal R}_1$ et ${\cal R}_2$ engendrent~$\cal R$ et agissent respectivement sur chacune des moiti\'es de ${\cal A}'$, on peut prolonger~${\cal A}'$ de mani\`ere $\cal R$-\'equivariante. On obtient alors un $\cal R$-champ de graphes bor\'elien connexes~$\cal A$ ; les sous-relations ${\cal R}_1$ et ${\cal R}_2$ \'etant en produit libre, ceci assure que les fibres de $\cal A$ soient des arbres. Autrement dit, $\cal A$ est un $\cal R$-arboretum tel que l'image de la section bor\'elienne $d$ soit un domaine fondamental de l'action de $\cal R$ sur~${\cal A}^0$, ce qui assure que ${\cal A}^0$ s'identifie au $\cal R$-espace fibr\'e standard canonique.
\findemo

\smallskip

\petiteremarque
Voici une deuxi\`eme d\'emonstration de la proposition pr\'ec\'edente, tr\`es similaire \`a la pr\'ec\'edente. Consid\'erons $Gr_1$ et $Gr_2$ des arborages de ${\cal R}_1$ et ${\cal R}_2$ respectivement et notons $Gr$ la r\'eunion de $Gr_1$ et $Gr_2$ : c'est une partie bor\'elienne de~$X \times X$, sym\'etrique, qui ne rencontre pas la diagonale, autrement dit un graphage. Puisque $Gr_1$ et $Gr_2$ engendrent respectivement ${\cal R}_1$ et ${\cal R}_2$ et que ces derni\`eres engendrent~$\cal R$, on en d\'eduit que $Gr$ engendre $\cal R$, autrement dit que $Gr$ est un graphage de $\cal R$. Enfin, s'il existait un cycle de longueur $n \geq 2$ dans le graphage $Gr$, on en d\'eduirait alors un $n$-uplet r\'eduit $\liste{x}$ d'\'el\'ements de $X$ tel que $x_n=x_1$, ce qui est exclu. Et ainsi, $Gr$ est un arborage de $\cal R$.

\bigskip

Consid\'erons une relation d'\'equivalence bor\'elienne $\cal R$ sur $X$ qui est le produit libre des sous-relations ${\cal R}_1$ et ${\cal R}_2$ d\'efinies sur $X$. Nous allons construire un $\cal R$-arboretum sur $X$ canoniquement associ\'e \`a cette d\'ecomposition. D\'esignons par $\tf{{\cal R}}{{\cal R}_1}$ et $\tf{{\cal R}}{{\cal R}_2}$ les $\cal R$-espaces fibr\'es standards canoniques associ\'es aux sous-relations ${\cal R}_1$ et ${\cal R}_2$. De plus, nous avons vu qu'il existait, \`a valeurs dans chacun de ces $\cal R$-espaces fibr\'es standards, un morphisme surjectif d\'efini sur le $\cal R$-espace fibr\'e standard canonique (cf. lem. \ref{morphisme fibr\'es} avec $A=X$, $s_1=d$ et $s_2=d_{{\cal R}_i}$), qui envoie la section diagonale $d$ sur la section bor\'elienne saturante $d_{{\cal R}_i}$. Notons $o$ le morphisme surjectif du $\cal R$-espace fibr\'e standard canonique dans $\tf{{\cal R}}{{\cal R}_1}$ et $t$ celui \`a valeurs dans $\tf{{\cal R}}{{\cal R}_2}$.

\medskip

Convenons d'appeler {\it espace des ar\^etes orient\'ees} le $\cal R$-espace fibr\'e standard canonique et {\it espace des sommets} la r\'eunion disjointe des $\cal R$-espaces fibr\'es standards~$\tf{{\cal R}}{{\cal R}_1}$ (espace des sommets de couleur 1) et $\tf{{\cal R}}{{\cal R}_2}$ (espace des sommets de couleur~2). Nous d\'esignerons \'egalement par {\it sommet origine} le morphisme de $\cal R$-espaces fibr\'es standards $o$ et par {\it sommet terminal} le morphisme $t$. Enfin, il n'est pas difficile de v\'erifier la compatibilit\'e des actions de $\cal R$ entre les espaces de sommets et d'ar\^etes. Ainsi nous avons d\'efini un $\cal R$-champ de graphes bor\'elien. La proposition qui suit assure qu'il s'agisse en fait d'un arboretum : cette construction est analogue \`a la construction de l'arbre de Bass-Serre associ\'e \`a un produit libre de deux groupes (cf. \cite{MR0476875}).

\begin{proposition}
Les fibres du $\cal R$-champ de graphes bor\'elien construit ci-dessus \`a partir d'une d\'ecomposition de $\cal R$ en produit libre de deux sous-relations d\'efinies sur~$X$ sont des arbres.
\end{proposition}
\label{arboretum canonique associe produit libre de sous-relations}
Dans la suite, nous dirons que le $\cal R$-arboretum construit pr\'ec\'edemment est le $\cal R$-arboretum canonique {\it associ\'e} \`a ${\cal R}={\cal R}_1 \star {\cal R}_2$.

\medskip

\demo
Il s'agit de v\'erifier que, pour tout \'el\'ement $x$ de $X$, la fibre de $x$ est un arbre, c'est-\`a-dire connexe et sans circuit (puisqu'on sait d\'ej\`a qu'elle est non vide). Soit $P$ et $Q$ deux sommets distincts dans le graphe au-dessus de $x$. Par d\'efinition de l'espace des sommets, $P$ et~$Q$ correspondent \`a des ${\cal R}_1$ ou ${\cal R}_2$-classes contenues dans la $\cal R$-classe de $x$. D\'esignons par $p$ et $q$ des repr\'esentants des classes correspondantes : ce sont donc des \'el\'ements $\cal R$-\'equivalents de $X$. Puisque~${\cal R}_1$ et ${\cal R}_2$ engendrent $\cal R$, il existe un $n$-uplet d'\'el\'ements que l'on peut supposer r\'eduit dans la $\cal R$-classe de $x$ dont le premier \'el\'ement est $p$ et le dernier $q$. Dans le graphe au-dessus de $x$, ce $n$-uplet correspond \`a un chemin joignant les sommets $P$ et $Q$.

\smallskip

Supposons par l'absurde qu'il existe un circuit dans le graphe au-dessus de $x$ dont on note $P_0$,... $P_{n-1}$ les sommets ($P_n=P_0$). Pour $i$ dans $\intn{0}{n-1}$, les sommets~$P_i$ et $P_{i+1}$ sont par hypoth\`ese reli\'es : il existe donc un \'el\'ement $x_{i+1}$ appartenant \`a~$X$ dont la ${\cal R}_1$-classe correspond \`a $P_i$ et la ${\cal R}_2$-classe \`a $P_{i+1}$ (ou l'inverse). De plus, par construction, deux sommets cons\'ecutifs sont de couleurs distinctes. Par suite, on construit un $n$-uplet $\liste{x}$ r\'eduit. Mais $x_1$ et $x_n$ \'etant ${\cal R}_j$-\'equivalents puisque les sommets $P_0$ et $P_n$ sont \'egaux, ceci est absurde.
\findemo

\smallskip

Par construction, le $\cal R$-arboretum canonique associ\'e \`a ${\cal R}={\cal R}_1 \star {\cal R}_2$ poss\`ede une section bor\'elienne d'ar\^etes dont l'image est un domaine fondamental pour l'action de $\cal R$ sur l'espace des ar\^etes orient\'ees et dont les sections bor\'eliennes de sommets associ\'ees ont leurs images qui forment une partition de l'espace des sommets en deux. De plus, les stabilisateurs de ces sections bor\'eliennes de sommets sont ${\cal R}_1$ et~${\cal R}_2$. Par suite, notons que pour toute ar\^ete de l'espace des ar\^etes, l'un de ses sommets est de couleur 1 et l'autre de couleur 2.

\bigskip

Nous allons d\'esormais d\'emontrer la r\'eciproque de ce r\'esultat et voir que les produits libres de sous-relations sont caract\'eris\'es par de telles actions : en un certain sens, la relation d'\'equivalence bor\'elienne engendr\'ee sur l'espace des ar\^etes est la plus simple possible.

\begin{theoreme}\label{d\'ecomposition des produits libres}
Soit $\cal R$ une relation d'\'equivalence bor\'elienne sur $X$. Supposons que~$\cal R$ agisse sur un arboretum $\cal A$ dont l'espace des ar\^etes orient\'ees ${\cal A}^{1+}$ est le $\cal R$-espace fibr\'e standard canonique et tel que les satur\'es des images des sections bor\'eliennes de sommets $o(d)$ et $t(d)$ forment une partition bor\'elienne de l'espace des sommets en deux. Alors $\cal R$ est le produit libre des stabilisateurs de $o(d)$ et de $t(d)$ :
\centers{${\cal R}=\text{Stab}_{\cal R}(o(d)) \star \text{Stab}_{\cal R}(t(d)).$}
\end{theoreme}

\smallskip

\demo
Notons ${\cal R}_1$ et ${\cal R}_2$ les stabilisateurs de $o(d)$ et $t(d)$ respectivement et commen\c{c}ons par montrer que $\cal R$ est engendr\'ee par ${\cal R}_1$ et ${\cal R}_2$. Soit $x$ et $y$ deux \'el\'ements $\cal R$-\'equivalents de~$X$ et consid\'erons l'image $(x,y) \cdot d(y)$ par le couple d'\'el\'ements $(x,y)$ de l'ar\^ete $d(y)$. Les ar\^etes $d(x)$ et $(x,y) \cdot d(y)$ appartiennent \`a la fibre~${\cal A}_x$ de $x$ et il existe donc une (unique) g\'eod\'esique $c$ entre les sommets $o(d(x))$ et $t((x,y) \cdot d(y))$ qui sont distincts par hypoth\`ese. Si $n$ est la longueur de ce chemin, d\'esignons par $\liste{w}$ les ar\^etes de $c$. Chacune d'entre elles appartient \`a l'orbite d'une unique ar\^ete appartenant \`a l'image de $d$ ou \`a l'image de $\overline{d}$ et, \'etant donn\'e deux ar\^etes cons\'ecutives, l'une appartient \`a l'orbite de l'image de $d$ et l'autre \`a l'orbite de l'image de $\overline{d}$. Puisque l'espace des ar\^etes orient\'ees est par hypoth\`ese le $\cal R$-espace fibr\'e standard canonique gauche, chacune des ar\^etes $w_i$ est un couple d'\'el\'ements $(x,x_i)$. On en d\'eduit ainsi un $n$-uplet $\liste{x}$ d'\'el\'ements de $X$ dont le premier co\"{i}ncide avec $x$, le dernier avec $y$ et tel que deux \'el\'ements successifs soient ${\cal R}_j$-\'equivalents pour un certain $j$.

\smallskip

Soit $x$ et $y$ deux \'el\'ements ${\cal R}_j$-\'equivalents de $X$. L'image du sommet $t(d(y))$ par le couple d'\'el\'ements $(x,y)$ est un sommet \`a distance unit\'e du sommet origine de~$d(x)$. Si $\liste{x}$ d\'esigne un $n$-uplet r\'eduit, on en d\'eduit par r\'ecurrence que le sommet terminal de $(x_1,x_n) \cdot d(x_n)$ est au moins \`a distance unit\'e du sommet origine de $d(x_1)$ : ceci exclut la possibilit\'e que $x_1$ soit \'egal \`a $x_n$.
\findemo

\bigskip

On d\'eduit du th\'eor\`eme pr\'ec\'edent le cas o\`u l'espace des ar\^etes orient\'ees est engendr\'e par l'image d'une section partielle dont le stabilisateur est trivial.

\begin{corollaire}
Soit $\cal A$ un $\cal R$-arboretum sur $X$. On suppose qu'il existe une section partielle $d$ (d\'efinie sur $A$) de l'espace des ar\^etes orient\'ees telle que $s(A)$ soit un domaine fondamental de ${\cal R}_{{\cal A}^{1+}}$. En outre, on suppose que les satur\'es des images des sections partielles de sommets $o(d)$ et $t(d)$ forment une partition bor\'elienne de l'espace des sommets en deux. Alors la restriction de $\cal R$ \`a~$A$ est le produit libre de ses sous-relations $\text{Stab}_{\cal R}(o(d))$ et $\text{Stab}_{\cal R}(t(d))$.
\findemo
\end{corollaire}

En particulier, on en d\'eduit qu'une telle relation d'\'equivalence bor\'elienne est stablement orbitalement \'equivalente \`a un produit libre de deux sous-relations.

\subsection{Produit amalgam\'e de deux sous-relations}

Nous allons g\'en\'eraliser au cas des produits amalgam\'es suivant une sous-relation commune (cf. d\'ef. \ref{d\'efinition produits amalgam\'es}) les r\'esultats obtenus dans le paragraphe pr\'ec\'edent. Puis nous terminons cette partie en donnant une application qui illustre que les m\'ethodes g\'eom\'etriques que nous d\'eveloppons permettent de d\'emontrer simplement des r\'esultats de nature alg\'ebrique.

\medskip

Soit $\cal R$ une relation d'\'equivalence bor\'elienne sur $X$ engendr\'ee par deux sous-relations ${\cal R}_1$ et~${\cal R}_2$. Soit ${\cal R}_3$ une sous-relation commune \`a ${\cal R}_1$ et~${\cal R}_2$. Soit un entier naturel $n \geq 2$. Un $n$-uplet $\liste{x}$ d'\'el\'ements de $X$ est dit {\it r\'eduit} si, pour tout~$i$ de $\intn{1}{n-1}$ :

\begin{itemize}
\item chaque $(x_i, x_{i+1})$ est un couple d'\'el\'ements ${\cal R}_1$-\'equivalents ou ${\cal R}_2$-\'equivalents et deux couples successifs ne sont pas ${\cal R}_j$-\'equivalents pour le m\^eme $j$ ;
\item aucun couple d'\'el\'ements $(x_i,x_{i+1})$ n'est ${\cal R}_3$-\'equivalent d\`es que $n>2$ ;
\item $x_1 \neq x_2 \qsiq n=2$.
\end{itemize}

\begin{definition}[Gaboriau, \cite{MR1728876}]\label{d\'efinition produits amalgam\'es}
La relation d'\'equivalence bor\'elienne $\cal R$ est le produit amalgam\'e des sous-relations ${\cal R}_1$ et ${\cal R}_2$ suivant la sous-relation ${\cal R}_3$ si, pour tout $n$-uplet r\'eduit $\liste{x}$ d'\'el\'ements de $X$, $x_n$ est diff\'erent de $x_1$. Dans ce cas, on note
\centers{${\cal R}={\cal R}_1 \star_{{\cal R}_3} {\cal R}_2.$}
\end{definition}

\petiteremarque
Dans le cas o\`u la sous-relation ${\cal R}_3$ est triviale, la notion de produit amalgam\'e des sous-relations ${\cal R}_1$ et~${\cal R}_2$ suivant la sous-relation ${\cal R}_3$ co\"{i}ncide avec celle de produit libre des sous-relations~${\cal R}_1$ et~${\cal R}_2$ (cf. d\'ef. \ref{d\'efinition produit libre}).

\smallskip

\exemple
Si un groupe d\'enombrable $\Gamma$ est le produit amalgam\'e de deux groupes~$\Gamma_1$ et $\Gamma_2$ suivant un sous-groupe commun $\Gamma_3$ et agit librement sur $X$, la relation d'\'equivalence bor\'elienne ${\cal R}_{\alpha}$ engendr\'ee par les orbites de l'action $\alpha$ de $\Gamma$ est le produit amalgam\'e des sous-relations ${\cal R}_{{\alpha}_1}$ et ${\cal R}_{{\alpha}_2}$ suivant la sous-relation ${\cal R}_{{\alpha}_3}$.

\medskip

\`A toute relation d'\'equivalence bor\'elienne $\cal R$ sur $X$ qui est le produit amalgam\'e des sous-relations ${\cal R}_1$ et ${\cal R}_2$ suivant la sous-relation ${\cal R}_3$, nous allons associer canoniquement un $\cal R$-arboretum ${\cal A}\pa{{\cal R}_1 \star_{{\cal R}_3} {\cal R}_2}$ sur $X$. D\'esignons par $\tf{{\cal R}}{{\cal R}_1}$, $\tf{{\cal R}}{{\cal R}_2}$ et $\tf{{\cal R}}{{\cal R}_3}$ les $\cal R$-espaces fibr\'es standards canoniques associ\'es aux paires $({\cal R},{\cal R}_1)$, $({\cal R},{\cal R}_2)$ et $({\cal R},{\cal R}_3)$. Or ${\cal R}_3$ \'etant une sous-relation commune de ${\cal R}_1$ et ${\cal R}_2$, on d\'eduit du lemme~\ref{morphisme fibr\'es} un morphisme surjectif~$o$ de $\tf{{\cal R}}{{\cal R}_3}$ dans $\tf{{\cal R}}{{\cal R}_1}$ envoyant la section bor\'elienne saturante~$d_{{\cal R}_3}$ sur la section bor\'elienne saturante $d_{{\cal R}_1}$, et de m\^eme un morphisme surjectif $t$ de~$\tf{{\cal R}}{{\cal R}_3}$ dans $\tf{{\cal R}}{{\cal R}_2}$ envoyant $d_{{\cal R}_3}$ sur $d_{{\cal R}_2}$.

\medskip

On d\'efinit l'espace des ar\^etes orient\'ees comme \'etant $\tf{{\cal R}}{{\cal R}_3}$ et l'espace des sommets comme la r\'eunion disjointe de $\tf{{\cal R}}{{\cal R}_1}$ (espace des sommets de couleur 1) et~$\tf{{\cal R}}{{\cal R}_2}$ (espace des sommets de couleur 2). L'application sommet origine est d\'efinie comme \'etant le morphisme de $\cal R$-espaces fibr\'es standards $o$ et l'application sommet terminal est d\'efinie comme \'etant le morphisme $t$. La compatibilit\'e des actions de $\cal R$ entre les espaces de sommets et d'ar\^etes assure que nous avons ainsi d\'efini un $\cal R$-champ de graphes bor\'elien $\cal A$. Nous allons maintenant d\'emontrer que les fibres de $\cal A$ sont des arbres, ce qui montrera que $\cal A$ est un $\cal R$-arboretum. Pour ceci, nous allons donner une caract\'erisation dynamique des produits amalgam\'es de sous-relations.

\begin{theoreme}\label{caract\'erisation des produits amalgam\'es}
Soit $\cal R$ une relation d'\'equivalence bor\'elienne et $\cal A$ un $\cal R$-champ de graphes bor\'elien sur $X$. On suppose que l'espace des ar\^etes orient\'ees ${\cal A}^{1+}$ est un espace fibr\'e standard homog\`ene dont on note $s : A \longrightarrow {\cal A}^{1+}$ une section partielle saturante. De plus, on suppose que les satur\'es des images des sections partielles de sommets $o(s)$ et $t(s)$ forment une partition $\cal R$-invariante de l'espace des sommets en deux. Les assertions suivantes sont alors \'equivalentes :
\begin{enumerate}
\item $\cal A$ est un $\cal R$-arboretum ;
\item $\rest{\cal R}{A}$ est le produit amalgam\'e des sous-relations $\text{Stab}_{\cal R}(o(s))$ et $\text{Stab}_{\cal R}(t(s))$ suivant la sous-relation $\text{Stab}_{\cal R}(s)$.
\end{enumerate}
\end{theoreme}

\demo
Puisque $A$ est un domaine complet de $\cal R$, quitte \`a consid\'erer le $\rest{\cal R}{A}$-arboretum sur $A$, on peut supposer que $A=X$ et se ramener au cas d'une section bor\'elienne saturante de l'espace des ar\^etes orient\'ees. Notons ${\cal R}_3=\text{Stab}_{\cal R}(s)$ le stabilisateur de $s$, et ${\cal R}_1=\text{Stab}_{\cal R}(o(s))$ et ${\cal R}_2=\text{Stab}_{\cal R}(t(s))$ les stabilisateurs des sections bor\'eliennes de sommets origines $o(s)$ et de sommets terminaux $t(s)$ correspondantes. Le th\'eor\`eme est alors une cons\'equence des deux lemmes suivants.

\begin{lemme}
Les fibres de $\cal A$ sont connexes si et seulement si $\cal R$ est engendr\'ee par~${\cal R}_1$ et ${\cal R}_2$.
\end{lemme}

\demo
Soit ${\cal A}'$ le champ de sous-graphes bor\'elien de $\cal A$ sur $X$ dont les fibres sont les composantes connexes de l'image de la section bor\'elienne d'ar\^etes $s$. Soit~${\cal R}'$ la sous-relation de~$\cal R$ qui laisse invariant le champ de sous-graphes bor\'elien ${\cal A}'$ : deux \'el\'ements $x$ et $y$ de $X$ sont ${\cal R}'$-\'equivalents si l'image par le couple d'\'el\'ements $(y,x)$ de~${\cal A}'_x$ est \'egale \`a ${\cal A}'_y$. Enfin on note ${\cal R}''$ la sous-relation de $\cal R$ engendr\'ee par~${\cal R}_1$ et ${\cal R}_2$. Si deux \'el\'ements $x$ et $y$ de $X$ sont ${\cal R}_1$ ou ${\cal R}_2$-\'equivalents, les ar\^etes $s(y)$ et $(y,x) \cdot s(x)$ ont un sommet en commun. L'ar\^ete $(y,x) \cdot s(x)$ appartient donc \`a ${\cal A}'_y$ et par suite
\centers{$(y,x) \cdot {\cal A}'_x={\cal A}'_y,$}
\noindent
ce qui prouve que ${\cal R}''$ est une sous-relation de ${\cal R}'$.

\smallskip

Notons ${\cal B}_1$ le champ de sous-graphes bor\'elien de $\cal A$ dont les ar\^etes sont les \'el\'ements du ${\cal R}''$-satur\'e de l'image de la section bor\'elienne d'ar\^etes $s$. Consid\'erons ${\cal B}_2$ le champ de sous-graphes bor\'elien de $\cal A$ dont les ar\^etes dans la fibre d'un \'el\'ement~$x$ de $X$ sont les $(x,y) \cdot s(y)$, o\`u $y$ est un \'el\'ement $\cal R$-\'equivalent \`a $x$ et n'appartenant pas \`a la ${\cal R}''$-classe de $x$. Comme l'image de $s$ est un domaine complet de la relation d'\'equivalence bor\'elienne ${\cal R}_{{\cal A}^{1+}}$ engendr\'ee sur l'espace des ar\^etes orient\'ees, la r\'eunion des champs de sous-graphes bor\'eliens disjoints ${\cal B}_1$ et ${\cal B}_2$ est \'egale \`a $\cal A$. Ceci implique que ${\cal B}_1$ contient ${\cal A}'$, et donc que ${\cal R}'$ est une sous-relation de ${\cal R}''$.

Les fibres de $\cal A$ sont connexes si et seulement si les champs de graphes bor\'elien~${\cal A}'$ et $\cal A$ sont \'egaux, c'est-\`a-dire si et seulement si $\cal R$ et ${\cal R}'$ sont \'egales, et d'apr\`es ce qui pr\'ec\`ede ${\cal R}'$ et ${\cal R}''$ co\"{i}ncident.
\findemo

\begin{lemme}
Pour qu'aucune fibre de $\cal A$ ne contienne de circuit, il faut et il suffit que pour tout $(n+1)$-uplet r\'eduit $\liste{x}$ d'\'el\'ements de $X$, l'\'el\'ement $x_n$ soit diff\'erent de $x_1$.
\end{lemme}

\demo
Pour tout \'el\'ement $x$ de $X$, la donn\'ee d'un chemin $c=\liste{c}$ ($n \geq 1$) sans aller-retour dans la fibre ${\cal A}_x$ de $x$ est \'equivalente \`a la donn\'ee d'un $n$-uplet $\liste{x}$ r\'eduit d'\'el\'ements de $X$ dans la $\cal R$-classe de $x$. En effet, chacune des ar\^etes de $c$ appartient \`a l'orbite d'une ar\^ete appartenant \`a l'image de $s$ ou \`a l'image de $\overline{s}$, donc pour tout $i$, il existe $x_i$ dans la $\cal R$-classe de $x$ tel que $c_i=(x,x_i) \cdot s(x_i)$ ou $c_i=(x,x_i) \cdot \overline{s}(x_i)$. De plus, \'etant donn\'e deux ar\^etes cons\'ecutives, l'une appartient \`a l'orbite de l'image de $s$ et l'autre \`a l'orbite de l'image de $\overline{s}$ d'apr\`es l'hypoth\`ese faite sur l'espace des sommets de $\cal A$. R\'eciproquement, un tel $n$-uplet d\'efinit bien un chemin sans aller-retour dans la fibre de $x$ dont les ar\^etes sont les $(x,x_i) \cdot s(x_i)$.

Enfin, le chemin $c$ est un circuit (c'est-\`a-dire tel que les sommets origine $o(c)$ et terminal $t(c)$ soient les m\^emes) si et seulement si les \'el\'ements $x_1$ et $x_n$ sont ${\cal R}_1$ ou ${\cal R}_2$-\'equivalents.
\findemo

\medskip

Comme cons\'equence du th\'eor\`eme pr\'ec\'edent et du fait d\'ej\`a vu que
\centers{$\text{Stab}_{\cal R}(\tilde{d_3})={\cal R}_3 \qquad \text{Stab}_{\cal R}(o(\tilde{d_3}))={\cal R}_1 \qqetqq \text{Stab}_{\cal R}(t(\tilde{d_3}))={\cal R}_2,$}
\noindent
nous en d\'eduisons le r\'esultat annonc\'e :

\begin{corollaire}\label{champ de graphes=arboretum pour produits amalgames}
Si une relation d'\'equivalence bor\'elienne $\cal R$ est le produit amalgam\'e de ses sous-relations ${\cal R}_1$ et ${\cal R}_2$ suivant la sous-relation ${\cal R}_3$, alors ${\cal A}\pa{{\cal R}_1 \star_{{\cal R}_3} {\cal R}_2}$ est un $\cal R$-arboretum.
\end{corollaire}

Nous dirons que $\cal A$ est le {\it $\cal R$-arboretum canonique associ\'e} \`a~${\cal R}={\cal R}_1 \star_{{\cal R}_3} {\cal R}_2$.

\bigskip

Pour terminer ce paragraphe, donnons une application de ce qui pr\'ec\`ede.

\begin{proposition}
Soit $\cal R$ une relation d'\'equivalence bor\'elienne qui est le produit amalgam\'e des sous-relations ${\cal R}_1$ et ${\cal R}_2$ suivant la sous-relation ${\cal R}_3$. Supposons que~$\cal S$ soit une sous-relation de $\cal R$ telle que, pour tout \'el\'ement $\phi$ de $[[\cal R]]$, les intersections des sous-relations $\phi^{-1} {\cal R}_1 \phi$ et $\phi^{-1} {\cal R}_2 \phi$ avec $\cal S$ soient lisses. Alors~$\cal S$ est une sous-relation arborable de $\cal R$.
\end{proposition}

\demo
 Consid\'erons le $\cal R$-arboretum canonique $\cal A$ associ\'e \`a la d\'ecomposition de $\cal R$ et montrons que l'action de $\cal S$ sur ${\cal A}^0$ est quasi-libre. Soit $s$ une section partielle de ${\cal A}^0$. Quitte \`a d\'ecouper en deux le domaine de d\'efinition de $s$, on peut supposer que $s$ est une section partielle de $\tf{{\cal R}}{{\cal R}_1}$ (ou $\tf{{\cal R}}{{\cal R}_2}$). Autrement dit, de mani\`ere g\'en\'erale, si $s$ est une section partielle de ${\cal A}^0$ d\'efinie sur $A \subset X$, alors $A$ est la r\'eunion disjointe des parties bor\'eliennes $A_1$ et $A_2$ telles que $s(A_1) \subset {\cal A}^{0,1}$ et $s(A_2) \subset {\cal A}^{0,2}$, o\`u~${\cal A}^{0,1}$ et ${\cal A}^{0,2}$ sont les espaces de sommets de $\cal A$ de couleurs 1 et 2 respectivement. La sous-relation $\cal S$ agit sur ${\cal A}^{0,1}=\tf{{\cal R}}{{\cal R}_1}$ et, par suite, le $\cal S$-stabilisateur de $s$ est l'intersection de~$\cal S$ avec une sous-relation de $\cal R$ (d\'efinie sur le domaine de d\'efinition de $s$) stablement conjugu\'ee \`a une restriction de ${\cal R}_1$ (cf. rem.~2 de la prop. \ref{stabilisateurs SOE}). Une telle intersection \'etant lisse par hypoth\`ese, on en d\'eduit que~$\cal S$ agit quasi-librement sur $\cal A$ et que $\cal S$ est arborable (cf. th. \ref{quasi-libre=arborable}).
\findemo

\section{$\cal R$-arboretums et d\'ecompositions de $\cal R$}
\seclabel{sec:cas g\'en\'eral}
\label{cas g\'en\'eral}

Dans cette partie, nous nous int\'eressons au cas g\'en\'eral d'une action d'une relation d'\'equivalence bor\'elienne $\cal R$ sur un arboretum $\cal A$ sur $X$. Nous sommes amen\'es \`a introduire la notion de {\it graphe de relations} (d\'ef. \ref{definition graphe de relations}) et les techniques que nous allons utiliser prolongent celles que nous avons d\'evelopp\'ees dans le cas des actions quasi-libres (cf. \S \ \ref{actions quasi-libres}). En particulier, nous \'etudions comment reconstruire $\cal R$ en termes de produits amalgam\'es \`a partir des stabilisateurs de certaines sections partielles de $\cal A$. Nous en d\'eduisons un th\'eor\`eme de d\'ecomposition (th. \ref{theoreme de Kurosh}) pour les sous-relations d'un produit libre d\'enombrable ${\cal R} = {\star}_{i \in I} {\cal R}_i$ de sous-relations ${\cal R}_i$ ($i \in I$), ainsi que la structure de toute restriction de $\cal R$ \`a une partie bor\'elienne $Y$ de $X$ (th. \ref{theoreme restriction}).

\subsection{Graphes de relations et d\'esingularisations}

Dans la d\'emonstration de la proposition \ref{arbre d'isomorphismes partiels} lors de notre \'etude des actions quasi-libres sur un arboretum, nous avons implicitement construit un {\it arbre d'isomorphismes partiels}, c'est-\`a-dire un arbre d\'enombrable \`a chaque sommet duquel est attach\'e un espace bor\'elien standard et \`a chaque ar\^ete un isomorphisme partiel entre des parties bor\'eliennes des espaces bor\'eliens standards port\'es par les extr\'emit\'es de cette ar\^ete. La d\'efinition suivante g\'en\'eralise cette id\'ee.

\begin{definition}[Graphe de relations]\label{definition graphe de relations}
Un graphe de relations $G^r$ est la donn\'ee d'un graphe d\'enombrable $G$, d'une relation d'\'equivalence bor\'elienne sur un espace bor\'elien standard pour chaque sommet et chaque ar\^ete de $G$ de sorte que les relations d'\'equivalence bor\'eliennes port\'ees par une ar\^ete et son ar\^ete oppos\'ee soient \'egales, ainsi que, pour chaque ar\^ete $a$ de $G$, d'un morphisme injectif de la relation d'\'equivalence bor\'elienne port\'ee par $a$ dans la relation d'\'equivalence bor\'elienne port\'ee par le sommet terminal de $a$. Si $G$ est un arbre, nous dirons que $G^r$ est un arbre de relations.
\end{definition}

\petiteremarque
Dans le cas particulier o\`u toutes les relations d'\'equivalence bor\'eliennes port\'ees par les sommets de $G$ sont triviales, nous dirons que $G^r$ est un {\it graphe d'isomorphismes partiels} (souvent not\'e $G^{ip}$) puisque la donn\'ee de $G^r$ permet d'attacher canoniquement \`a chaque ar\^ete $a$ de $G$ un isomorphisme partiel $\phi_a$ dont la source et le but sont des parties bor\'eliennes des espaces bor\'eliens standards port\'es par les sommets origine et terminal de $a$, de sorte que $\phi_a$ et $\phi_{\bar{a}}$ soient inverses l'un de l'autre. Si de plus, tous les espaces bor\'eliens standards port\'es par les sommets de $G$ sont des singletons, alors la donn\'ee de $G^{ip}$ est \'equivalente \`a celle du graphe $G$.

\medskip
\label{relation arbre d'isomorphismes partiels}
\'Etant donn\'e un arbre d'isomorphismes partiels $G^{ip}$, consid\'erons la r\'eunion disjointe des espaces bor\'eliens standards port\'es par les sommets de $G$. Les isomorphismes partiels canoniquement port\'es par les ar\^etes de $G$ engendrent une relation d'\'equivalence bor\'elienne lisse $\cal S$ sur cette r\'eunion disjointe. Un graphage $Gr$ de $\cal S$ est d\'efini de la mani\`ere suivante : un couple d'\'el\'ements $(x,y)$ appartient \`a $Gr$ si $x$ et~$y$ appartiennent aux espaces bor\'eliens standards port\'es par des sommets adjacents de~$G$ et si $y$ est l'image de $x$ par l'isomorphisme partiel port\'e par l'une des deux ar\^etes (inverses l'une de l'autre) joignant ces deux sommets. De plus, \'etant donn\'e une \'enum\'eration des sommets de $G$, il n'est pas difficile d'exhiber par r\'ecurrence un domaine fondamental de $\cal S$. Notons~$D_1$ l'espace borélien standard port\'e par le premier sommet : puisque le graphe $G$ sous-jacent \`a $G^{ip}$ est un arbre, $D_1$ ne rencontre au plus qu'une seule fois chaque classe de $\cal S$. Pour tout $k \geq 1$, ayant construit une partie bor\'elienne $D_k$ de la r\'eunion des espaces bor\'eliens standards port\'es par les $k$ premiers sommets de $G$ qui ne rencontre au plus qu'une seule fois chaque classe de $\cal S$, $D_{k+1}$ est par d\'efinition la r\'eunion de $D_k$ et du compl\'ementaire du satur\'e de $D_k$ dans l'espace bor\'elien standard port\'e par le sommet $k+1$ de $G$. La r\'eunion croissante des $D_k$ est un domaine fondamental de cette relation.

\smallskip

\label{arbre d'iso. part. enracine}
Un arbre d'isomorphismes partiels $G^{ip}$ est dit {\it enracin\'e} s'il est muni d'un sommet, appel\'e {\it racine}, dont l'espace bor\'elien standard attach\'e soit un domaine fondamental de cette relation d'\'equivalence bor\'elienne. Notons que dans ce cas, l'isomorphisme partiel port\'e par une ar\^ete dont le sommet origine est celui de ses deux sommets le plus proche de la racine est surjectif. Un tel arbre d'isomorphismes partiels enracin\'e d\'efinit naturellement un arboretum sur l'espace bor\'elien standard~$Y$ port\'e par sa racine. Si $\cal A$ est un champ de graphes bor\'elien sur $X$ et $G^{ip}$ un arbre d'isomorphismes partiels enracin\'e, une {\it repr\'esentation} de $G^{ip}$ dans $\cal A$ est la donn\'ee d'un isomorphisme d'espaces bor\'eliens standards entre $Y$ et une partie bor\'elienne~$X_Y$ de~$X$ et d'un isomorphisme $\chi$ entre l'arboretum d\'efini par $G^{ip}$ sur sa racine $Y$ et un sous-arboretum de la restriction de $\cal A$ \`a $X_Y$. On dit alors que $\chi$ est une repr\'esentation de $G^{ip}$ dans~$\cal A$ et que les images (par l'isomorphisme de repr\'esentation $\chi$) des espaces bor\'eliens standards port\'es par les sommets et les ar\^etes de $G^{ip}$ sont des {\it repr\'esentations} dans~$\cal A$ de ces espaces bor\'eliens standards.

\medskip

\'Etant donn\'e une relation d'\'equivalence bor\'elienne $\cal R$ sur $X$, consid\'erons d\'esormais le cas d'un $\cal R$-champ de graphes bor\'elien $\cal A$ sur $X$. Si $\chi$ est une repr\'esentation d'un arbre d'isomorphismes partiels enracin\'e $G^{ip}$ dans $\cal A$, l'action de $\cal R$ sur $\cal A$ donne naissance \`a des relations d'\'equivalence bor\'eliennes sur les repr\'esentations dans $\cal A$ des espaces bor\'eliens standards de sommets et d'ar\^etes de $G^{ip}$ : ce sont respectivement les restrictions de ${\cal R}_{{\cal A}^0}$ et ${\cal R}_{{\cal A}^1}$. Une repr\'esentation de $G^{ip}$ dans $\cal A$ induit ainsi une structure d'arbre de relations $G^r$ sur $G^{ip}$.

\smallskip

Via le foncteur d'{\it oubli}, \`a tout arbre de relations $G^r$ est canoniquement associ\'e un arbre d'isomorphismes partiels $G^{ip}$. Nous dirons qu'un arbre de relations $G^r$ est {\it enracin\'e} si l'arbre d'isomorphismes partiels $G^{ip}$ sous-jacent \`a $G^r$ l'est. Une {\it repr\'esentation d'un arbre de relations enracin\'e $G^r$ dans $\cal A$} est une repr\'esentation $\chi$ de $G^{ip}$ dans $\cal A$ telle que $\chi$ induise des isomorphismes de relations d'\'equivalence bor\'eliennes entre les relations d'\'equivalence bor\'eliennes port\'ees par les sommets et les ar\^etes de $G^r$ et les restrictions de ${\cal R}_{{\cal A}^0}$ et ${\cal R}_{{\cal A}^1}$ aux repr\'esentations dans $\cal A$ des espaces bor\'eliens standards port\'es par les sommets et les ar\^etes de $G^{ip}$.

\begin{definition}[Arboretum de repr\'esentants]\label{d\'efinition arboretum de repr\'esentants}
Soit $\cal A$ un $\cal R$-champ de graphes bor\'elien sur $X$. Un arboretum de repr\'esentants de l'action de $\cal R$ sur $\cal A$ est la donn\'ee d'un arbre de relations enracin\'e $G^r$ et d'une repr\'esentation $\chi$ de $G^r$ dans $\cal A$ telle que les  ${\cal R}_{{\cal A}^0}$-satur\'es des repr\'esentations des espaces bor\'eliens standards port\'es par les sommets de $G^r$ forment une partition de l'espace des sommets~${\cal A}^0$.
\end{definition}

\petiteremarque
\'Etant donn\'e un relev\'e d'un arbre maximal de l'espace quotient d'une action d'un groupe sur un graphe connexe, on d\'efinit naturellement un arbre de groupes en consid\'erant les stabilisateurs des sommets et des ar\^etes de ce relev\'e (cf.~\cite{MR0476875}). Un arboretum de repr\'esentants correspond \`a une donn\'ee analogue dans le cas des relations d'\'equivalence bor\'eliennes.

\medskip

Comme dans le cas d'une action quasi-libre sur un arboretum (cf. prop.~\ref{arbre d'isomorphismes partiels}), on d\'emontre l'existence d'un arboretum de repr\'esentants pour toute action de~$\cal R$ sur un arboretum $\cal A$ sur $X$.

\begin{proposition}\label{arboretum de repr\'esentants}
Pour tout $\cal R$-arboretum $\cal A$ sur $X$, il existe un arboretum de repr\'esentants de l'action de $\cal R$ sur $\cal A$.
\end{proposition}

\demo
L'id\'ee de la preuve de l'existence d'un arboretum de repr\'esentants est essentiellement la m\^eme que dans le cas des actions quasi-libres (cf. prop. \ref{arbre d'isomorphismes partiels}). Dans la d\'emonstration de la proposition \ref{arbre d'isomorphismes partiels}, nous prenions soin \`a chaque \'etape de consid\'erer des sections partielles de sommets (et donc d'ar\^etes) dont les stabilisateurs \'etaient triviaux, ce qui \'etait possible car ${\cal R}_{{\cal A}^0}$ et ${\cal R}_{{\cal A}^1}$ \'etaient des relations d'\'equivalence bor\'eliennes lisses par hypoth\`ese, pour ainsi construire un arbre d'isomorphismes partiels. Dans le cas g\'en\'eral o\`u ${\cal R}_{{\cal A}^0}$ et ${\cal R}_{{\cal A}^1}$ ne sont plus suppos\'ees lisses, il n'est plus possible {\it a priori} d'assurer que les stabilisateurs des sections partielles choisies soient triviaux \`a chaque \'etape de la construction. C'est ainsi que nous sommes contraints de passer des arbres d'isomorphismes partiels aux arbres de relations : le stabilisateur d'une section partielle de sommets ajout\'ee \`a l'\'etape $n+1$ et le stabilisateur de la section partielle d'ar\^etes correspondante d\'efinissent des relations d'\'equivalence bor\'eliennes port\'ees par un sommet terminal attach\'e par une ar\^ete (et son ar\^ete oppos\'ee) \`a l'arbre de relations construit \`a l'\'etape $n$.
\findemo

\smallskip

\petiteremarque
La premi\`ere section de sommets construite \'etant d\'efinie sur un domaine complet $A$ quelconque de $\cal R$, on peut choisir $A=X$ (c'est-\`a-dire une section bor\'elienne $s$ d\'efinie sur tout $X$) et ainsi supposer que l'espace bor\'elien standard port\'ee par la racine s'identifie (via $s$) \`a $X$, ce que nous ferons par la suite.

\bigskip

Nous allons introduire une derni\`ere notion, celle de {\it d\'esingularisation d'une action}. \'Etant donn\'e une action d'une relation d'\'equivalence bor\'elienne $\cal R$ sur un arboretum $\cal A$ sur $X$, nous avons vu (prop. \ref{arboretum de repr\'esentants}) l'existence d'un arboretum de repr\'esentants $G^r$ de cette action mais cette notion n' \og encode \fg pas toute l'information de l'action de $\cal R$ sur $\cal A$ dans le cas g\'en\'eral puisque les  ${\cal R}_{{\cal A}^0}$-satur\'es des repr\'esentations des espaces bor\'eliens standards port\'es par les ar\^etes de $G^r$ ne recouvrent pas {\it a priori} l'espace des ar\^etes ${\cal A}^1$.

\begin{definition}[D\'esingularisation d'une action]\label{definition desingularisation d'une action}
Une d\'esingularisation de l'action de $\cal R$ sur $\cal A$ est la donn\'ee d'un graphe orient\'e de relations $G^r$ (c'est-\`a-dire de graphe sous-jacent $G$ orient\'e), d'un sous-arbre maximal enracin\'e $A_m$ de $G$ et d'un isomorphisme partiel $\phi_a$ de~$[[\cal R]]$ pour chaque ar\^ete $a$ de l'orientation de $G$ n'appartenant pas \`a $A_m$ tel que :
\begin{itemize}
\item l'arbre de relations enracin\'e $A^r_m$ induit par $G^r$ sur $A_m$ soit un arboretum de repr\'esentants de l'action de $\cal R$ sur $\cal A$ tel que l'espace bor\'elien standard port\'e par la racine de $A_m$ s'identifie \`a $X$ (notons $\chi$ l'isomorphisme de repr\'esentation) ;
\item pour toute ar\^ete $a$ de l'orientation de $G$ n'appartenant pas \`a $A_m$, il existe une section partielle d'ar\^etes $s_a$ d\'efinie sur la source de $\phi_a$, dont les sommets origines appartiennent \`a la repr\'esentation dans $\cal A$ de l'espace bor\'elien standard port\'e par le sommet origine de $a$ et telle que les images par $\phi_a$ des sommets terminaux appartiennent \`a la repr\'esentation dans $\cal A$ de l'espace bor\'elien standard port\'e le sommet terminal de $a$. Comme pour les sommets, nous dirons que l'image de cette section partielle est une repr\'esentation de l'espace bor\'elien standard port\'e par l'ar\^ete ;
\item pour toute ar\^ete $a$ de l'orientation de $G$ n'appartenant pas \`a $A_m$, $\chi$ induit un isomorphisme entre $\text{Stab}_{\cal R}(s_a)$ identifi\'e \`a une sous-relation de ${\cal R}_{\chi(o(a))}$ et l'image par le morphisme injectif associ\'e \`a $\bar{a}$ de la relation d'\'equivalence bor\'elienne port\'ee par $\bar{a}$ dans la relation d'\'equivalence bor\'elienne port\'ee par $t(\bar{a})$, ainsi qu'un isomorphisme entre $\phi_a \text{Stab}_{\cal R}(s_a) \phi_a^{-1}$ identifi\'e \`a une sous-relation de ${\cal R}_{\chi(t(a))}$ et l'image par le morphisme injectif associ\'e \`a~$a$ de la relation d'\'equivalence bor\'elienne port\'ee par $a$ dans la relation d'\'equivalence bor\'elienne port\'ee par $t(a)$ ;
\item les ${\cal R}_{{\cal A}^1}$-satur\'es des repr\'esentations des espaces bor\'eliens standards port\'es par les ar\^etes de $G^r$ forment une partition de l'espace des ar\^etes ${\cal A}^1$.
\end{itemize}
\end{definition}

D\'emontrons \`a pr\'esent l'existence d'une d\'esingularisation pour toute action d'une relation d'\'equivalence bor\'elienne sur un arboretum. Ce r\'esultat d'existence est analogue \`a la construction d'un graphe de groupes associ\'e \`a une action d'un groupe sur un arbre.

\begin{theoreme}\label{desingularisation d'une action}
Soit $\cal A$ un $\cal R$-arboretum sur $X$. Alors il existe une d\'esingularisation de l'action de~$\cal R$ sur~$\cal A$.
\end{theoreme}

\noindent
{\it Remarque 1 :}
Comme nous allons le voir dans la construction ci-dessous, pour toute ar\^ete $a$ de l'orientation de $G$ n'appartenant pas \`a $A_m$, le graphe de l'isomorphisme partiel $\phi_a$ ne rencontre pas la diagonale de $X \times X$.

\smallskip

\noindent
{\it Remarque 2 :}
Tout comme dans le cas d'une action d'un groupe sur un arbre, il n'y a pas de d\'esingularisation canonique et nous sommes amen\'es \`a faire de nombreux choix lors de notre construction. D'ailleurs il n'y a d\'ej\`a pas de choix canonique d'un graphe sous-jacent \`a une d\'esingularisation, contrairement au cas classique o\`u le graphe sous-jacent est l'espace quotient de l'action sans inversion du groupe sur l'arbre.

\medskip

\demo
La proposition \ref{arboretum de repr\'esentants} assure l'existence d'un arboretum de repr\'esentants, c'est-\`a-dire un arbre de relations enracin\'e $A^r_m$ et d'un isomorphisme de repr\'esentation $\chi$ de $A^r_m$ dans $\cal A$, tel que l'espace bor\'elien standard port\'e par la racine de $A_m$ s'identifie \`a $X$. Ainsi les orbites sous l'action de~$\cal R$ des repr\'esentations des espaces bor\'eliens standards port\'es par les sommets de $A^r_m$ forment une partition $\cal R$-invariante de l'espace des sommets ${\cal A}^0$. D\'esignons par $C$ le compl\'ementaire dans~${\cal A}^1$ des $\cal R$-orbites des repr\'esentations des espaces bor\'eliens standards port\'es par les ar\^etes de $A^r_m$. Si $C$ est vide, alors $A^r_m$ est d\'ej\`a une d\'esingularisation de l'action de~$\cal R$ sur~$\cal A$. Sinon, la partie bor\'elienne~$C'$ de $C$ constitu\'ee des ar\^etes dont le sommet origine appartient \`a l'arboretum de repr\'esentants est non vide car les fibres de $\cal A$ sont des graphes connexes et que $\cal R$ agit en pr\'eservant les distances sur ${\cal A}^0$. Nous allons construire une famille de sections partielles de $C$ dont les satur\'es des images forment une partition bor\'elienne de $C$.

Soit $P$ un sommet de $A^r_m$ et $s$ la section partielle de sommets de $\cal A$ repr\'esentant l'espace bor\'elien standard port\'e par $P$. Pour tout sommet $Q$ de $A^r_m$, consid\'erons la partie bor\'elienne $C_{(P,Q)}$ de $C$ constitu\'ee des ar\^etes dont le sommet origine appartient \`a l'image de $s$ et dont le sommet terminal appartient au satur\'e de la repr\'esentation de l'espace bor\'elien standard port\'e par $Q$. Si $C_{(P,Q)}$ est non vide, c'est un espace fibr\'e standard sur $\pi(C_{(P,Q)})$ que nous pouvons exhauster d'apr\`es la remarque suivant le lemme \ref{partition} \`a l'aide de sections partielles dont les satur\'es des images forment une partition bor\'elienne de $C_{(P,Q)}$. Chacune de ces sections partielles d'ar\^etes donne lieu \`a une ar\^ete dans $A^r_m$ entre les sommets $P$ et $Q$. La relation d'\'equivalence bor\'elienne~${\cal R}_a$ induite par ${\cal R}_{{\cal A}^1}$ sur l'image de l'une de ces sections partielles~$a$ de~${\cal A}^1$ s'identifie naturellement \`a une sous-relation de la repr\'esentation de la relation d'\'equivalence bor\'elienne port\'ee par $P$ et, quitte \`a consid\'erer une restriction de $a$ \`a un domaine complet de son stabilisateur, ${\cal R}_a$ s'injecte comme sous-relation de la relation d'\'equivalence bor\'elienne port\'ee par $Q$ via un isomorphisme partiel de $[[\cal R]]$ qui envoie la section de sommets terminaux $t(a)$ dans la repr\'esentation de l'arboretum de repr\'esentants.
\findemo

\bigskip

\'Etant donn\'e une d\'esingularisation $(G^r,A_m,\chi)$ de l'action d'une relation d'\'equivalence bor\'elienne $\cal R$ sur un arboretum $\cal A$ sur $X$, nous allons d\'esormais d\'emontrer quelques propri\'et\'es de cette d\'esingularisation. Rappelons que la repr\'esentation de l'arboretum de repr\'esentants de l'action de~$\cal R$ sur $\cal A$ est un sous-arboretum de $\cal A$ sur $X$ et que l'isomorphisme partiel port\'e par une ar\^ete de $A^{ip}_m$ dont le sommet origine est celui de ses deux sommets le plus proche de la racine est surjectif. Enfin, rappelons que les relations d'\'equivalence bor\'eliennes port\'ees par les sommets de $G^r$ s'identifient via $\chi$ \`a des sous-relations de $\cal R$ et qu'\`a chaque ar\^ete $a$ de l'orientation de $G$ n'appartenant pas \`a~$A_m$ correspond un isomorphisme partiel $\phi_a$ de $[[\cal R]]$. On d\'efinit $\phi_{\overline{a}}=\phi_a^{-1}$ pour toute ar\^ete $a$ de l'orientation de $G \setminus A_m$ et, pour toute ar\^ete~$a$ de $A_m$, $\phi_a$ comme la restriction de l'identit\'e de $X$ sur la projection dans $X$ de la repr\'esentation de l'espace bor\'elien standard port\'e par le sommet origine de $a$.

\smallskip

Soit $P$ et $Q$ deux sommets adjacents dans l'arbre maximal $A_m$ et supposons que~$P$ soit le plus proche des deux du sommet racine. Identifi\'ees via $\chi$  \`a des sous-relations de $\cal R$ sur les parties bor\'eliennes $X_P \supset X_Q$ de $X$, notons ${\cal R}_P$ et ${\cal R}_Q$ les relations d'\'equivalence bor\'eliennes port\'ees par les sommets $P$ et $Q$. De m\^eme, identifi\'ee \`a une sous-relation de $\cal R$ d\'efinie sur $X_Q$, d\'esignons par ${\cal R}_{P,Q}$ la relation d'\'equivalence bor\'elienne port\'ee par l'ar\^ete de $A_m$ joignant $P$ \`a $Q$ .

\begin{lemme}
La relation d'\'equivalence bor\'elienne engendr\'ee par $\rest{{\cal R}_P}{X_Q}$ et ${\cal R}_Q$ sur~$X_Q$ est le produit amalgam\'e de $\rest{{\cal R}_P}{X_Q}$ et ${\cal R}_Q$ suivant la sous-relation ${\cal R}_{P,Q}$.
\end{lemme}

\demo
Consid\'erons dans $\cal A$ les repr\'esentations des espaces bor\'eliens standards port\'es par les sommets~$P$ et $Q$, ainsi que la repr\'esentation de l'espace bor\'elien standard port\'e par l'ar\^ete de~$A^r_m$ joignant $P$ \`a $Q$ ; cette derni\`ere d\'efinit sur $X_Q \subset X$ une section partielle d'ar\^etes $s$. Notons~$\cal S$ la relation d'\'equivalence bor\'elienne engendr\'ee par $\rest{{\cal R}_P}{X_Q}$ et ${\cal R}_Q$ sur $X_Q$ et consid\'erons l'orbite de l'image de $s$ sous l'action de~$\cal S$. Par d\'efinition de $\cal S$, on obtient pour tout \'el\'ement $x$ de $X_Q$ un sous-arbre~${\cal A}'_x$ de~${\cal A}_x$. Ainsi, $\cal S$ est une relation d'\'equivalence bor\'elienne qui agit sur l'arboretum~${\cal A}'$ sur~$X_Q$, lequel poss\`ede une section bor\'elienne d'ar\^etes dont l'image est un domaine complet de ${\cal S}_{{\cal A}'^1}$ et tel que les satur\'es des sections bor\'eliennes de sommets origines~$o(s)$ et terminaux $t(s)$ forment une partition bi-color\'ee de l'espace des sommets~${\cal A}^{'0}$. Le r\'esultat d\'ecoule alors du th\'eor\`eme \ref{caract\'erisation des produits amalgam\'es}.
\findemo

\medskip

Le lemme pr\'ec\'edent s'\'etend sans difficult\'e au cas de deux sommets quelconques distincts. Pour ceci, consid\'erons la g\'eod\'esique dans $A^r_m$ joignant les deux sommets~$P$ et $Q$ et supposons que l'intersection $X_{P,Q} \subset X$ des espaces bor\'eliens standards port\'es par ces sommets ne soit pas vide. D\'esignons par ${\cal R}_{P,Q}$ l'intersection de toutes les relations d'\'equivalence bor\'eliennes port\'ees par les ar\^etes de la g\'eod\'esique joignant~$P$ et~$Q$.

\begin{proposition}\label{produits amalgames}
La relation d'\'equivalence bor\'elienne $\cal S$ engendr\'ee par $\rest{{\cal R}_P}{X_{P,Q}}$ et $\rest{{\cal R}_Q}{X_{P,Q}}$ sur $X_{P,Q}$ est le produit amalgam\'e de $\rest{{\cal R}_P}{X_{P,Q}}$ et $\rest{{\cal R}_Q}{X_{P,Q}}$ suivant la sous-relation ${\cal R}_{P,Q}$.
\end{proposition}

\demo
Nous allons proc\'eder comme dans la d\'emonstration du lemme pr\'ec\'edent : consid\'erons les repr\'esentations dans $\cal A$ des espaces bor\'eliens standards port\'es par les sommets et les ar\^etes de la g\'eod\'esique joignant $P$ \`a $Q$ dans $A_m$. Partant du sous-arboretum ${\cal A}''$ de la restriction de $\cal A$ \`a $X_{P,Q}$, consid\`erons l'image de~${\cal A}''$ sous l'action de la relation d'\'equivalence bor\'elienne $\cal S$. On obtient ainsi un $\cal S$-arboretum~${\cal A}'$ sur $X_{P,Q}$. La contraction de ${\cal A}'$ suivant ${\cal A}''$ dans chaque fibre de $X_{P,Q}$ est une op\'eration bor\'elienne qui se prolonge par $\cal S$-\'equivariance : on obtient alors un nouvel $\cal S$-arboretum sur $X_{P,Q}$ satisfaisant les conditions du th\'eor\`eme~\ref{caract\'erisation des produits amalgam\'es} et dont le stabilisateur de la section d'ar\^etes ainsi construite n'est autre que l'intersection de toutes les relations d'\'equivalence bor\'eliennes port\'ees par les ar\^etes de la g\'eod\'esique joignant~$P$ \`a~$Q$.
\findemo

\medskip

Identifi\'ees via $\chi$ \`a des sous-relations de $\cal R$, les relations d'\'equivalence bor\'eliennes port\'ees par les sommets de $G$ engendrent une sous-relation ${\cal R}'$ de $\cal R$ d\'efinie sur $X$. Notons ${\cal R}''$ la sous-relation de $\cal R$ ({\it a priori} d\'efinie sur une partie bor\'elienne de $X$) engendr\'ee par les isomorphismes partiels~$\phi_a$ port\'es par chacune des ar\^etes $a$ de $G$. Nous allons d\'esormais voir que $\cal R$ est engendr\'ee par~${\cal R}'$ et~${\cal R}''$.

\begin{lemme}\label{intersection}
Soit $x$ et $y$ deux \'el\'ements $\cal R$-\'equivalents de $X$.
\begin{enumerate}
\item $(y,x) \cdot {\cal A}'_x \cap {\cal A}'_y \neq \emptyset$ si et seulement si $x$ et $y$ sont ${\cal R}_P$-\'equivalents pour la relation d'\'equivalence bor\'elienne ${\cal R}_P$ port\'ee par un sommet $P$ de $G$ ;
\item le ${\cal R}'$-satur\'e de ${\cal A}'$ est un ${\cal R}'$-arboretum ;
\item $(y,x) \cdot \pac{{\cal R}' \cdot {\cal A}'_x} \cap \pac{{\cal R}' \cdot {\cal A}'_y} \neq \emptyset$ si et seulement si $x$ et $y$ sont ${\cal R}'$-\'equivalents. Dans ce cas on a en fait $(y,x) \cdot \pac{{\cal R}' \cdot {\cal A}'_x} = \pac{{\cal R}' \cdot {\cal A}'_y}$.
\end{enumerate}
\end{lemme}

\demo
Le point 1 est clair par d\'efinition et le point 2 est une cons\'equence directe de la connexit\'e dans chaque fibre d\'eduite du point 1. On en d\'eduit alors directement le point 3.
\findemo

\begin{proposition}\label{relation arborable}
\begin{enumerate}
\item Toute relation d'\'equivalence bor\'elienne $R_P$ port\'ee par un sommet $P$ de $G$ est d'intersection triviale avec ${\cal R}''$ ;
\item ${\cal R}''$ est arborable.
\end{enumerate}
\end{proposition}

\demo
Le point 1 est une cons\'equence du fait que les fibres de $\cal A$ sont des arbres et donc ne contiennent pas de circuit de longueur 1. On d\'eduit du point 1 et du point 1 du lemme~\ref{intersection} que les images par des couples d'\'el\'ements ${\cal R}''$-\'equivalents de~${\cal A}'$ dans une fibre donn\'ee sont disjoints. Consid\'erons dans chaque fibre l'enveloppe convexe de tous ces translat\'es, autrement dit, l'enveloppe convexe du ${\cal R}''$-satur\'e de~${\cal A}'$. La contraction par ${\cal R}''$-\'equivariance suivant ${\cal A}'$ donne un nouvel arboretum sur lequel ${\cal R}''$ agit : le stabilisateur de la section diagonale apr\`es contraction est trivial. En effet, si $x$ et $y$ appartiennent au stabilisateur de la section diagonale, alors $(y,x) \cdot {\cal A}'_x$ est \'egal \`a ${\cal A}'_y$ et on a d\'ej\`a vu que ceci implique $x=y$. On obtient finalement un ${\cal R}''$-arboretum dont une section bor\'elienne de sommets d\'efinie sur $X$ est un domaine fondamental pour l'action sur l'espace des sommets : il s'agit donc d'une action lisse, ce qui prouve l'arborabilit\'e de ${\cal R}''$.
\findemo

\begin{lemme}\label{lemme section sommets saturante}
Soit $\cal S$ une relation d'\'equivalence bor\'elienne agissant sur un arboretum $\cal B$ sur un espace bor\'elien standard $Y$. Supposons qu'il existe une section bor\'elienne $d$ de sommets saturante de stabilisateur ${\cal S}'$. Si $(s_a)_{a \in \underline{A}}$ est une famille d\'enombrable de sections partielles d'ar\^etes, dites {\it extra-ar\^etes}, telle que :
\begin{enumerate}
\item les sommets origines de ces extra-ar\^etes appartiennent \`a $d(Y)$ ;
\item pour chacune de ces sections d'ar\^etes $s_a$ d\'efinie sur $A_a$, il existe un isomorphisme partiel $\phi : A_a \longrightarrow B_a$ de $\cal S$ tel que l'image des sommets terminaux de $s_a(A_a)$ par $\phi$ soit exactement $d(B_a)$ ;
\item la r\'eunion de ces extra-ar\^etes soit un domaine complet de la relation d'\'equivalence bor\'elienne sur l'espace des ar\^etes,
\end{enumerate}
\noindent
alors $\cal S$ est engendr\'ee par ${\cal S}'$ et ${\cal S}''$, o\`u ${\cal S}''$ d\'esigne la relation d'\'equivalence bor\'elienne engendr\'ee par le L-graphage $(\phi_a)_{a \in \underline{A}}$.
\end{lemme}

\demo
Soit $x$ et $y$ deux \'el\'ements $\cal S$-\'equivalents de $Y$. Int\'eressons-nous \`a la g\'eod\'esique (de longueur $k$) joignant $d(x)$ \`a l'image par $(x,y)$ de $d(y)$, et notons $a_1$ la premi\`ere ar\^ete de cette g\'eod\'esique de sommet origine $d(x)$. Si $a_1$ est une extra-ar\^ete, alors il existe un \'el\'ement $x_1$ de $Y$ qui soit ${\cal S}''$-\'equivalent \`a $x$ et tel que l'orbite du sommet terminal de $a_1$ contienne $d(x_1)$. Puisque~$\cal S$ agit sur $\cal B$ en pr\'eservant les distances, la longueur de la g\'eod\'esique joignant $d(x_1)$ \`a l'image par $(x_1,y)$ de $d(y)$ est de longueur $k-1$. Sinon, quitte \`a consid\'erer un \'el\'ement ${\cal S}'$-\'equivalent \`a $x$, on peut supposer que $a_1$ est l'image d'une extra-ar\^ete par un des isomorphismes partiels~$\phi_a$. Dans ce cas, il existe $x_1$ tel que l'image de $d(x_1)$ par $\phi_a$ soit exactement $d(x)$ : l\`a encore, la longueur de la g\'eod\'esique joignant $d(x_1)$ \`a l'image par $(x_1,y)$ de $d(y)$ est alors de longueur $k-1$. Soit $n$ un entier $\geq 1$. \'Etant donn\'e $0 \leq k \leq n-1$ et une suite d'\'el\'ements $\cal S$-\'equivalents $\liste[k]{x}$ \`a $x$, l'argument pr\'ec\'edent permet de conclure que si la longueur de la g\'eod\'esique joignant $d(x_{n-k})$ \`a l'image par $(x_{n-k},y)$ de $d(y)$ est de longueur $n-k$, alors, quitte \`a consid\'erer d'abord un \'el\'ement ${\cal S}'$-\'equivalent, il existe un \'el\'ement ${\cal S}''$-\'equivalent $x_{n-k-1}$ \`a $x_{n-k}$ tel que la longueur de la g\'eod\'esique $d(x_{n-k-1})$ \`a l'image par $(x_{n-k-1},y)$ de $d(y)$ soit de longueur $n-k-1$. D'o\`u le r\'esultat.
\findemo

\begin{proposition}\label{generateurs}
La relation d'\'equivalence bor\'elienne $\cal R$ est engendr\'ee par~${\cal R}'$ et~${\cal R}''$.
\end{proposition}
\demo
Consid\`ereons la contraction de $\cal A$ suivant le ${\cal R}'$-satur\'e de ${\cal A}'$. D'apr\`es le point~2 du lemme \ref{intersection}, on obtient un $\cal R$-arboretum satisfaisant aux hypoth\`eses du lemme pr\'ec\'edent, le stabilisateur de la section bor\'elienne de sommets $d$ \'etant pr\'ecis\'ement ${\cal R}'$ d'apr\`es le point 3 du m\^eme lemme.
\findemo

\subsection{Sous-relations d'un produit libre}

Nous allons utiliser les r\'esultats pr\'ec\'edents pour \'etudier les sous-relations d'une relation d'\'equivalence bor\'elienne qui admet une d\'ecomposition en produit libre d'un ensemble d\'enombrable de sous-relations. En th\'eorie des groupes, le th\'eor\`eme de Kurosh (\cite{MR0071422}) assure qu'un sous-groupe~$H$ du produit libre $G$ d'une famille de groupes $(G_z)_{z \in Z}$ est isomorphe au produit libre de son intersection avec des conjugu\'es des~$G_z$ convenablement index\'es et d'un sous-groupe libre de $G$. Nous obtenons un th\'eor\`eme de d\'ecomposition (th. \ref{theoreme de Kurosh}) dans le contexte des relations d'\'equivalence bor\'eliennes pour les sous-relations d'un produit libre ${\cal R} = {\star}_{i \in I} {\cal R}_i$ d'un ensemble d\'enombrable de sous-relations ${\cal R}_i$, ainsi que la structure d'une restriction de $\cal R$ \`a toute partie bor\'elienne $Y$ de $X$ (th. \ref{theoreme restriction}).

\begin{definition}[Produit amalgam\'e d\'enombrable]\label{definition produit amalgame denombrable}
La relation d'\'equivalence bor\'elienne $\cal R$ sur $X$ admet une d\'ecomposition en produit amalgam\'e d\'enombrable des sous-relations ${\cal R}_i$ ($i \in I$) suivant la sous-relation commune ${\cal R}'$ si $\cal R$ est engendr\'ee par la famille d\'enombrable des ${\cal R}_i$ et si, pour tout $n$-uplet $\liste{x}$ d'\'el\'ements de $X$ tels que $x_n=x_1$ et $(x_k,x_{k+1}) \in {\cal R}_{i_k}$, il existe $1 \leq k \leq {n-1}$ tel que $(x_k,x_{k+1}) \in {\cal R}'$. Dans ce cas, on note
\centers{${\cal R}=\star_{{\cal R}'} {\cal R}_i.$}
\end{definition}

En particulier, si ${\cal R}'$ est triviale dans la d\'efinition pr\'ec\'edente, on obtient la notion de produit libre d\'enombrable de sous-relations.

\begin{theoreme}\label{theoreme de Kurosh}
Soit $\cal R$ une relation d'\'equivalence bor\'elienne sur $X$, produit libre d\'enombrable de sous-relations ${\cal R}_i$ ($i \in I$). Alors
\centers{${\cal S} = {\star}_{i \in I}\pa{{\star}_{k_i \in K(i)} {\cal S}_{k_i}} \star {\cal T},$}
\noindent
o\`u, pour tout $k_i$ d'un ensemble d\'enombrable $K(i)$, il existe un \'el\'ement $\phi_{k_i}$ de $[[\cal R]]$ d\'efini sur une partie bor\'elienne $A_{k_i}$ de $X$ tel que
\centers{${\cal S}_{k_i} = \pa{\phi_{k_i}^{-1} \rest{{\cal R}_{i}}{\phi_{k_i}(A_{k_i})} \phi_{k_i}} \cap {\cal S},$}
\noindent
et o\`u $\cal T$ est une sous-relation arborable de $\cal S$. De plus, pour tout $i$ de $I$, il existe $k_i$ dans~$K(i)$ tels que
\centers{$A_{k_i}=X \qqetqq {\cal S}_{k_i}={\cal R}_i \cap {\cal S}.$}
\end{theoreme}

\petiteremarque
Le r\'esultat pr\'ec\'edent s'\'etend sans difficult\'e suppl\'ementaire au cas d'un produit amalgam\'e d\'enombrable de sous-relations ${\cal R}_i$ suivant une sous-relation commune~${\cal R}'$ sous l'hypoth\`ese que l'intersection de $\cal S$ avec toute conjugu\'ee de~${\cal R}'$ soit lisse.

\smallskip

\demo
Commen\c{c}ons par d\'eduire le cas g\'en\'eral du cas particulier d'un produit libre de deux sous-relations ${\cal R}_1$ et ${\cal R}_2$. Donnons-nous un ordre sur $I$ (index\'e par les entiers naturels non nuls, \'eventuellement un nombre fini si $I$ est fini) et commen\c{c}ons par isoler le premier facteur dans la d\'ecomposition de $\cal R$ : on \'ecrit ainsi ${\cal R}={{\cal R}_1} \star_{{\cal R}'} {{\cal R}'_2}$, o\`u ${\cal R}'_2$ est le produit amalgam\'e des ${\cal R}_i$ restantes suivant ${\cal R}'$. On en d\'eduit que $\cal S$ est le produit libre de son intersection avec des $\cal R$-conjugu\'ees de restrictions de ${\cal R}_1$ et ${\cal R}'_2$ et d'une sous-relation arborable. On est ainsi amen\'e \`a s'int\'eresser \`a certaines sous-relations de ${\cal R}'_2={{\cal R}_2} \star_{{\cal R}'} {{\cal R}'_3}$ (dont ${\cal S} \cap {\cal R}'_2$), o\`u ${\cal R}_2$ est le deuxi\`eme facteur dans la d\'ecomposition initiale de $\cal R$ et ${\cal R}'_3$ le produit amalgam\'e des ${\cal R}_i$ restantes suivant~${\cal R}'$. On applique \`a nouveau le th\'eor\`eme dans le cas de deux facteurs \`a cette sous-relation, ce qui donne une d\'ecomposition faisant intervenir une sous-relation arborable et des sous-relations de $\cal R$-conjugu\'ees de restrictions de ${\cal R}_2$ et ${\cal R}'_3$. On continue ainsi la d\'ecomposition des sous-relations de $\cal R$-conjugu\'ees de restrictions des ${{\cal R}'_i}$ apparaissant et, en utilisant un argument diagonal et le fait que le produit libre de sous-relations arborables est encore arborable (cf. prop. \ref{produit libre arborable}, le cas d'une infinit\'e d\'enombrable de facteurs \'etant analogue), on en d\'eduit l'existence de ${\cal S}_{k_i}$ et de $\cal T$ comme souhait\'ees telles que ${\star}_{i \in I}\pa{{\star}_{k_i \in K(i)} {\cal S}_{k_i}} \star {\cal T}$ soit une sous-relation de $\cal S$. Enfin si deux \'el\'ements~$x$ et~$y$ de $X$ sont $\cal S$-\'equivalents, il existe un nombre fini de ${\cal R}_i$ tel que $(x,y)$ appartienne au produit libre d\'enombrable de ces ${\cal R}_i$, et donc $(x,y)$ appartient \`a ${\star}_{i \in I}\pa{{\star}_{k_i \in K(i)} {\cal S}_{k_i}} \star~{\cal T}$.

\medskip

Nous allons maintenant d\'emontrer le cas particulier d'un produit libre de deux sous-relations ${\cal R}_1$ et ${\cal R}_2$. Pour cela, consid\'erons le $\cal R$-arboretum canonique $\cal A$ associ\'e \`a la d\'ecomposition en produit libre de ${\cal R}={\cal R}_1 \star {\cal R}_2$ (cf. p. \pageref{arboretum canonique associe produit libre de sous-relations}). Pr\'ecisons une d\'esingularisation de l'action de $\cal S$ (d\'eduite de celle de $\cal R$) sur $\cal A$ en suivant la construction d\'evelopp\'ee dans la d\'emonstration de la proposition~\ref{arboretum de repr\'esentants} et en utilisant la partition naturelle de l'espace des sommets de $\cal A$ en deux parties. Autrement dit, \`a chaque \'etape de la construction d'un arboretum de repr\'esentants (dont on note~${\cal A}'$ sa repr\'esentation dans $\cal A$), les sections partielles de sommets que l'on consid\`ere sont des sections partielles de~$\tf{{\cal R}}{{\cal R}_1}$ ou de $\tf{{\cal R}}{{\cal R}_2}$ et la premi\`ere section bor\'elienne d'ar\^etes que l'on choisit est la diagonale de l'espace fibr\'e standard canonique des ar\^etes orient\'ees. Or nous avons d\'ej\`a vu (cf. rem. 2 de la prop. \ref{stabilisateurs SOE}) que le stabilisateur d'une section partielle de $\tf{{\cal R}}{{\cal R}_i}$ est stablement conjugu\'e \`a la restriction de ${\cal R}_i$ \`a une partie bor\'elienne de $X$. Quitte \`a restreindre son domaine de d\'efinition, on peut donc supposer \`a chaque \'etape que le stabilisateur de la section partielle de $\tf{{\cal R}}{{\cal R}_i}$ construite soit conjugu\'e \`a une restriction de~${\cal R}_i$ et que le stabilisateur de la section partielle d'ar\^etes ainsi d\'efini soit trivial : ceci est possible car l'action de $\cal R$ est lisse sur l'espace des ar\^etes. Enfin, puisque la premi\`ere section bor\'elienne d'ar\^etes consid\'er\'ee est la diagonale de~$\cal R$, on en d\'eduit que les deux premi\`eres sections de sommets sont les diagonales $d_{{\cal R}_1}$ et $d_{{\cal R}_2}$ de $\tf{{\cal R}}{{\cal R}_1}$ et $\tf{{\cal R}}{{\cal R}_2}$. Notons que ce qui pr\'ec\`ede s'\'etend sans difficult\'e suppl\'ementaire au cas o\`u~$\cal R$ est le produit amalgam\'e de ${\cal R}_1$ et ${\cal R}_2$ suivant la sous-relation commune ${\cal R}_3$ sous l'hypoth\`ese que l'intersection de $\cal S$ avec toute conjugu\'ee de ${\cal R}_3$ soit lisse. En effet cette derni\`ere hypoth\`ese assure que $\cal S$ agisse quasi-librement sur le $\cal R$-espace fibr\'e standard $\tf{{\cal R}}{{\cal R}_3}$.

\medskip

Contractons ${\cal A}'$ par $\cal R$-\'equivariance : on se ram\`ene ainsi au cadre du lemme \ref{lemme section sommets saturante} sous l'hypoth\`ese 3') plus forte qui assure que la r\'eunion des extra-ar\^etes soit un domaine fondamental de la relation d'\'equivalence bor\'elienne sur l'espace des ar\^etes, ce qui est une cons\'equence du fait que l'action de~$\cal R$ soit libre. Utilisant les notations du lemme \ref{lemme section sommets saturante} et compte tenu des propositions \ref{produits amalgames}, \ref{relation arborable} et \ref{generateurs}, il ne reste plus qu'\`a voir que ${\cal S}'$ et ${\cal S}''$ sont des sous-relations en produit libre de $\cal S$.

Commen\c{c}ons par remarquer que ${\cal S}'$ et ${\cal S}''$ sont d'intersection triviale. En effet, si deux \'el\'ements~$x$ et~$y$ sont ${\cal S}'$-\'equivalents, alors le sommet $d(x)$ a pour image le sommet $d(y)$ sous l'action de $(y,x)$ ; or si $x$ et $y$ sont ${\cal S}''$-\'equivalents, par d\'efinition, un des sommets \`a distance 1 de $d(x)$ est envoy\'e sur $d(y)$. Si $x \neq y$ \'etait $\pa{{\cal S}' \cap {\cal S}''}$-\'equivalents, on en d\'eduirait alors un cycle de longueur 1 attach\'e en~$d(y)$.

\smallskip

\petiteremarque Rappelons \`a pr\'esent que si on op\`ere par un couple d'\'el\'ements $(t,\phi_a(t))$ o\`u $t$ appartient \`a la source $A_a$ de $\phi_a$, cela signifie qu'au-dessus de $t$, il existe une ar\^ete du domaine fondamental dont le sommet origine est $d(t)$ et que le sommet terminal de cette ar\^ete est envoy\'e par $(\phi_a(t),t)$ sur $d(\phi(t))$.

\smallskip

Soit $v$ un sommet dans la fibre de $x$. Notons $a_1$ la premi\`ere ar\^ete (de sommet origine $d(x)$) de la g\'eod\'esique reliant $d(x)$ \`a $v$. Analysons les diff\'erentes situations possibles :

\begin{enumerate}
\item si on op\`ere par un couple d'\'el\'ements $(x,y)$ de ${\cal S}'$, alors la distance de $v$ \`a $d(x)$ dans la fibre de $x$ est pr\'eserv\'ee par l'action de $(y,x)$ ;

\item si on op\`ere par un couple d'\'el\'ements $(x,\phi_a(x))$ pour un certain $\phi_a$ :
\begin{itemize}
\item si $a_1$ est l'extra-ar\^ete qui d\'efinit $\phi_a$ (et donc appartient au domaine fondamental), alors la distance augmente de $-1$ ;
\item si $a_1$ n'est pas l'extra-ar\^ete qui d\'efinit $\phi_a$, alors la distance augmente de $+1$ ;
\end{itemize}

\item si on op\`ere par un couple d'\'el\'ements $(x=\phi_a(u),u)$ pour un certain $\phi_a$ :
\begin{itemize}
\item si l'ar\^ete oppos\'ee de $a_1$ est l'image de l'extra-ar\^ete qui d\'efinit $\phi_a$, alors la distance augmente de $-1$ ;
\item si l'ar\^ete oppos\'ee de $a_1$ n'est pas l'image de l'extra-ar\^ete qui d\'efinit $\phi_a$, alors la distance augmente de $+1$.
\end{itemize}
\end{enumerate}

Supposons un mouvement de la forme : $+1$ suivi de $-1$. L'\'el\'ement $(\phi_a(x),x) \cdot d(x)$ est un \'el\'ement \`a distance $+1$ de $d(\phi_a(x))$. Consid\'erons un mouvement de la forme~$-1$ avec un certain~$\psi_b$ ; d'apr\`es ce qui pr\'ec\`ede, il s'agit des situations 2-point 1 et 3-point~1. La situation 2-point 1 est impossible car $a_1$, qui est l'ar\^ete joignant $d(\phi(x))$ \`a $(\phi(x),x) \cdot d(x)$ est par d\'efinition telle que son oppos\'ee est l'image d'une ar\^ete du domaine fondamental : en particulier, ce n'est pas une extra-ar\^ete pouvant d\'efinir~$\psi_b$. Donc 3-point 1 est la seule possibilit\'e et l'unicit\'e due au domaine fondamental assure alors $\psi_b$ soit l'inverse de $\phi_a$. Ainsi tout mouvement de la forme $+1$ suivi de~$-1$ est trivial, c'est-\`a-dire de la forme : $x \longrightarrow \phi_a(x) \longrightarrow x$ (conclusion 1).

Supposons un mouvement de la forme : $+1$ suivi de $0$ suivi de $-1$. Les m\^emes arguments que ci-dessus prouvent que cela est impossible, si le $0$ correspond \`a deux points distincts dans ${\cal S}'$ (conclusion 2).

Donnons-nous une suite de points $\liste{x}$ avec $x_n = x_1$, successivement~${\cal S}'$ ou~${\cal S}''$-\'equivalents et \'etudions la distance \`a la section $d$ des images successives du sommet $d(x_0)$ sous l'action de ces couples d'\'el\'ements. Appliquons les transformations successives de $x_1$ \`a $x_n = x_1$, que l'on suppose \^etre une suite r\'eduite. La fonction distance est une fonction \`a valeurs dans $\Z$ qui \`a chaque \'etape avance de $0$, d'un entier positif ou d'un entier n\'egatif, qui part de $0$ et qui revient en~$0$ (les fibres de~$\cal B$ sont des arbres). Ce qui pr\'ec\`ede (cf. conclusion 1) assure qu'on se ram\`ene alors \`a un mouvement de la forme
\centers{$0 \longrightarrow (+1+1+1+1) \longrightarrow 0 \longrightarrow (+1+1+1) \longrightarrow 0 \longrightarrow (-1-1) \longrightarrow 0 \longrightarrow \ \text{etc}. $}
\noindent
l'uplet \'etant toujours r\'eduit si on regroupe chaque parenth\`ese. Mais nous avons vu (cf. conclusion~2) que la distance ne fait que cro\^{i}tre ou d\'ecro\^{i}tre strictement selon le signe de la premi\`ere valeur, ce qui est absurde car $x_n = x_1$.
\findemo

\medskip

Donnons \`a pr\'esent une deuxi\`eme application des r\'esultats pr\'ec\'edents concernant la restriction de ${\cal R}={\cal R}_1 \star {\cal R}_2$ \`a une partie bor\'elienne $A$ de $X$ qui est un domaine complet de ${\cal R}_1$ et ${\cal R}_2$. On retrouve alors les r\'esultats de Ioana-Peterson-Popa (\cite{ioana-peterson-popa}, prop. 7.4.2) obtenus dans le cadre mesur\'e.

\begin{proposition}
Soit ${\cal R}={\cal R}_1 \star {\cal R}_2$ une relation d'\'equivalence bor\'elienne sur $X$, produit libre des sous-relations ${\cal R}_1$ et ${\cal R}_2$. Si $A$ d\'esigne un domaine complet commun de ${\cal R}_1$ et ${\cal R}_2$, alors
\centers{$\rest{\cal R}{A}=\rest{{\cal R}_1}{A} \star \rest{{\cal R}_2}{A} \star {\cal T}$}
\noindent
o\`u $\cal T$ est une sous-relation arborable de $\cal R$.
\end{proposition}

\demo
Continuons de d\'esigner par $\cal A$ le $\cal R$-arboretum canonique associ\'e \`a la d\'ecomposition de $\cal R$. Consid\'erons l'action de $\rest{\cal R}{A}$ sur la restriction de $\cal A$ \`a~$A$. Par hypoth\`ese, l'image de $A$ par $d_{{\cal R}_i}$ est un domaine complet de ${\cal R}_{\tf{{\cal R}}{{\cal R}_i}}$ et, par cons\'equent, la r\'eunion de $d_{{\cal R}_1}(A)$ et de $d_{{\cal R}_2}(A)$ est un domaine complet de la restriction de~${\cal R}_{{\cal A}^0}$ \`a l'espace des sommets de $\rest{\cal A}{A}$. La construction d\'evelopp\'ee dans la d\'emonstration de la proposition \ref{arboretum de repr\'esentants} donne une d\'esingularisation de l'action de~$\rest{\cal R}{A}$ sur~$\rest{\cal A}{A}$ telle que le graphe de relations associ\'e ait exactement deux sommets portant les relations d'\'equivalence bor\'eliennes $\rest{{\cal R}_1}{A}$ et $\rest{{\cal R}_2}{A}$ et tel que les relations d'\'equivalence bor\'eliennes associ\'ees aux ar\^etes soient triviales. Les m\^emes arguments que dans la d\'emonstration du th\'eor\`eme~\ref{theoreme de Kurosh} permettent de conclure.
\findemo

\medskip

Formulons le cas g\'en\'eral (sans hypoth\`ese sur $A$) pour un produit libre d\'enombrable de sous-relations. Puisque $\rest{\cal R}{A}$ est une sous-relation de $\cal R$, on peut en particulier appliquer le th\'eor\`eme pr\'ec\'edent. Mais nous pouvons \^etre plus pr\'ecis encore.

\begin{theoreme}\label{theoreme restriction}
Soit $\cal R$ une relation d'\'equivalence bor\'elienne sur $X$, produit libre d\'enombrable de sous-relations ${\cal R}_i$ ($i \in I$) d\'efinies sur des parties bor\'eliennes $A_i$ de~$X$ telles que la r\'eunion des~$A_i$ soit \'egale \`a $X$. Pour toute partie bor\'elienne $A$ de $X$, la restriction $\rest{\cal R}{A}$ de $\cal R$ \`a $A$ admet une d\'ecomposition en produit libre d\'enombrable de sous-relations de la forme
\centers{$\rest{\cal R}{A} = {\star}_{i \in I}\pa{{\star}_{k_i \in K(i)} {\cal S}_{k_i}} \star {\cal T},$}
\noindent
o\`u, pour tout $k_i$ de l'ensemble d\'enombrable $K(i)$, il existe un \'el\'ement $\phi_{k_i}$ de $[[\cal R]]$ d\'efini sur une partie bor\'elienne $A_{k_i}$ de $X$ tel que
\centers{${\cal S}_{k_i} = \phi_{k_i}^{-1} \rest{{\cal R}_{i}}{\phi_{k_i}(A_{k_i})} \phi_{k_i},$}
\noindent
et o\`u $\cal T$ d\'esigne une sous-relation arborable de $\cal R$. De plus, pour tout $i$ de $I$, la r\'eunion des ${\cal R}_i$-satur\'es des buts $\phi_{k_i}(A_{k_i})$ des $\phi_{k_i}$ forme une partition (d\'enombrable et bor\'elienne) du domaine de d\'efinition $A_i$ de ${\cal R}_i$. Enfin, pour tout $i$ de $I$ tel que $A_i \cap A$ soit non vide, il existe $k_i$ dans~$K(i)$ tels que
\centers{$A_{k_i}=A_i \cap A \qqetqq {\cal S}_{k_i}=\rest{{\cal R}_i}{A_i \cap A}.$}
\end{theoreme}

\demo
Le cas g\'en\'eral se d\'eduit du cas de deux facteurs comme dans la d\'emonstration du th\'eor\`eme \ref{theoreme de Kurosh}. D\'esignant toujours par $\cal A$ le $\cal R$-arboretum canonique associ\'e \`a la d\'ecomposition de ${\cal R}={\cal R}_1 \star {\cal R}_2$, consid\'erons le $\rest{\cal R}{A}$-arboretum~$\rest{\cal A}{A}$ et d\'esingularisons l'action de $\rest{\cal R}{A}$ exactement comme dans la d\'emonstration du th\'eor\`eme \ref{theoreme de Kurosh} en utilisant la partition naturelle de l'espace des sommets de $\cal A$ en deux parties ; la premi\`ere section bor\'elienne d'ar\^etes que l'on choisit est bien entendu la restriction \`a $A$ de la diagonale de l'espace fibr\'e standard canonique des ar\^etes orient\'ees, de sorte que les deux premi\`eres sections de sommets soient les restrictions $\rest{d_{{\cal R}_1}}{A}$ et $\rest{d_{{\cal R}_2}}{A}$ \`a $A$ des diagonales $d_{{\cal R}_1}$ et $d_{{\cal R}_2}$ de $\tf{{\cal R}}{{\cal R}_1}$ et~$\tf{{\cal R}}{{\cal R}_2}$.

\smallskip

Le seul point qu'il reste \`a voir est que, pour tout $i$ de $I$, la r\'eunion des ${\cal R}_i$-satur\'es des $\phi_{k_i}(A_{k_i})$ forme une partition de $A_i$. Pour cela, rappelons que la r\'eunion des images $d_{{\cal R}_1}(X)$ et $d_{{\cal R}_2}(X)$ des diagonales~$d_{{\cal R}_1}$ et~$d_{{\cal R}_2}$ est un domaine complet de ${\cal R}_{{\cal A}^0}$. Par suite, quitte \`a consid\'erer sa restriction \`a un domaine complet de son $\cal R$-stabilisateur, pour toute section partielle~$s$ de sommets monochromes (de couleur $l=1, 2$) de $\rest{\cal A}{A}$ d\'efinie sur une partie bor\'elienne~$A_s$ de $A$, il existe un isomorphisme partiel $\phi_s : A_s \longrightarrow B_s$ de $[[{\cal R}]]$ tel que l'image sous l'action de $\phi_s$ de $s(A_s)$ soit exactement $d_{{\cal R}_l}(B_s)$. Si $s$ et $s'$ d\'esignent deux telles sections de sommets de m\^eme couleur $l$ dans la construction d'une d\'esingularisation de l'action de $\rest{\cal R}{A}$ (cf. d\'em. prop. \ref{arboretum de repr\'esentants}) , on en d\'eduit que les ${\cal R}_{{\cal A}^0}$-satur\'es de $d_{{\cal R}_l}(B_s)$ et $d_{{\cal R}_l}(B_{s'})$ sont disjoints, et par suite, que les ${\cal R}_l$-satur\'es de~$B_s$ et~$B_{s'}$ sont disjoints car~${\cal R}_l$ est pr\'ecis\'ement le $\cal R$-stabilisateur de $d_{{\cal R}_l}$.
\findemo

\providecommand{\bysame}{\leavevmode ---\ }
\providecommand{\og}{``}
\providecommand{\fg}{''}
\providecommand{\smfandname}{et}
\providecommand{\smfedsname}{\'eds.}
\providecommand{\smfedname}{\'ed.}
\providecommand{\smfmastersthesisname}{M\'emoire}
\providecommand{\smfphdthesisname}{Th\`ese}

\end{document}